\documentclass[10pt]{amsart}
\usepackage{array,amsmath,amssymb,amsthm,amsxtra,amscd,latexsym,graphics,psfrag,
epsfig}

\theoremstyle{plain}
\newtheorem{thm}{Theorem}
\newtheorem{cor}{Corollary}

\newtheorem{lemma}{Lemma}
\newtheorem{prop}{Proposition}

\newtheorem{property}{Property}
\newtheorem{lemma*}{Lemma}

\theoremstyle{definition}
\newtheorem{definition}{Definition}

\theoremstyle{remark}
\newtheorem*{notation}{Notation}
\newtheorem*{rmk}{Remark}

\numberwithin{equation}{section}

\newcommand{\gs}[2][]{GS_{(#1,#2)}} 
\newcommand{\vanish}[1]{}

\newcommand{\OO}{\mathcal{O}}
\newcommand{\PP}{\mathcal P}
\newcommand{\LL}{\mathcal L}
\newcommand{\UU}{\mathcal U}
\newcommand{\ncl}{\sim_l}

\begin{document}

\title{Partially ordered groups which act on oriented order trees}

\date{\today}


\author{Matthew Horak$^1$}
\address{Department of Mathematics\\
         Trinity College\\
         Hartford, CT 06106}
\email{Matthew.Horak@trincoll.edu}
\author{Melanie Stein$^2$}
\address{Department of Mathematics\\
         Trinity College\\
         Hartford, CT 06106}
\email{mstein@trincoll.edu}

\keywords{partially ordered group, 3-manifold, tree, order tree, 1-manifold, group action} \subjclass{20E08,20F60,06F15,20F65,57M60}

\begin{abstract}
It is well known that a countable group admits a left-invariant total order if and only if it acts faithfully on $\mathbb{R}$ by orientation preserving homeomorphisms. Such group actions are special cases of group actions on simply connected 1-manifolds, or equivalently, actions on oriented order trees. We characterize a class of left-invariant partial orders on groups which yield such actions, and show conversely that groups acting on oriented order trees by order preserving homeomorphism admit such partial orders as long as there is an action with a point whose stabilizer is left-orderable.

\end{abstract}

\maketitle

\section{Introduction}\label{S:intro}

It is a well known result that a countable group admits a left (or right)-invariant total order if and only if it acts faithfully by orientation preserving homeomorphisms on $\mathbb{R}$ (see Theorem 6.8 of \cite{Gh} for a proof). Many groups that arise naturally in topology are left-orderable, for example Artin groups, certain mapping class groups of Riemann surfaces, and many 3-manifold groups.  Boyer, Rolfsen and Wiest establish in \cite{BRW} that there are compact connected manifolds modelled on each of the eight 3-dimensional geometries with both orderable and nonorderable fundamental groups.  The first examples of non-orderable hyperbolic 3-manifold groups are given in \cite{RSS}, and the non-orderability of these groups is established by showing that they cannot act via faithful orientation preserving homeomorphisms on $\mathbb{R}$.  Now, $\mathbb{R}$ is a special case of a simply connected 1-manifold, and in fact the paper shows that these groups cannot act nontrivially on any oriented simply connected 1-manifold. Furthermore, the points of oriented simply-connected 1-manifolds have a natural order, just as the points of $\mathbb{R}$ do. In contrast to the total order of the points on the real line, the points of these more general manifolds are in general only partially ordered (see \cite{RSS}). The goal of this paper is to characterize the left invariant partial orders on groups that correspond to group actions on oriented simply-connected 1-manifolds, and their generalizations, oriented order trees, and to prove the analogue of the classical theorem relating left orderability of groups and faithful group actions on $\mathbb{R}$.

The paper is organized as follows. In section 2, we recall some of the relevant background information on simply connected 1-manifolds and their generalizations, order trees. We also extend some basic ideas about ends from Hausdorff trees to the more general setting of order trees. In section 3, we define the abstract partial orders we will be concerned with, which we call partial orders with rectifiable simply connected extensions (see Definitions ~\ref{simplyconnected}~(simply connected posets),~\ref{D:extension}~(extensions), and~\ref{D:finiteswitches}~(rectifiable)). We establish certain properties of sets equipped with such partial orders, which allow us to define a notion of ``betweenness" (see Definition~\ref{D:betweenness}) for arbitrary simply connected posets that agrees with the natural notion of ``between" for 1-manifolds. We also show that these orders naturally pass to both subgroups and quotient groups (see Theorem~\ref{T:quotientgetsorder}). In section 4 we use the partial order and its betweenness relation to show by construction that groups that admit left invariant partial orders with rectifiable simply connected extensions act on simply connected 1-manifolds by orientation preserving homeomorphisms (see Theorem~\ref{T:order.to.action}).  In section 5, we show in Theorems \ref{T:arbitrarytree} and \ref{T:weakstab} that groups which act on oriented order trees by orientation preserving homeomorphisms admit rectifiable left-invariant partial orders with simply connected extensions, as long as there is some point of the tree whose stabilizer is totally ordered.  In fact, groups that act minimally on simply connected 1-manifolds admit simply connected partial orders (see Theorems~\ref{T:groupgetsorder} and \ref{T:nontrivialstab}), whereas a group that acts minimally only on a more general oriented order tree may not admit such a partial order, though it does admit an order with formal extension.  

The study of treelike structures is common throughout the literature, in for example the theories of $\mathbb R$-trees, $\Lambda$-trees, protrees, dendrons, pretrees, etc.  A certainly non-exhaustive survey of the literature developing the theories of such structures includes,~\cite{AB}, \cite{Ba}, \cite{bow}, \cite{bc}, \cite{CM}, \cite{dG94}, \cite{p96}, and~\cite{sh1}.  In particular, the notion of ``betweenness" and the associated betweenness relation is central to the theory of pretrees. Bowditch and Crisp~\cite{bow},~\cite{bc}  have successfully used betweenness relations to generalize results for topological actions on $\mathbb R$-trees to actions on pretrees.  After we completed the constructions of Section 4, it came to our attention that the constructions of Bowditch and Crisp are very similar to our constructions.  Indeed, our betweenness relation on a set with a simply connected partial order does satisfy Bowditch's notion of betweenness, and hence a simply connected poset has the structure of a pretree. However, the constructions in~\cite{bc} and ~\cite{bow} applied to one of our groups would give an action on an $\mathbb R$-tree without orientation, since the pretrees themselves have no notion of order. One could give a proof of Theorem ~\ref{T:order.to.action} by strengthening the constructions of Bowditch and Crisp in the special case of a  simply connected partial order to carry the full poset information (referring the reader to those constructions), and then orienting the resulting $\mathbb R$-trees and proving that their actions preserve the orientations. However, including the extra structure on the poset from the beginning allows for a more straightforward construction. Hence, for reasons of both completeness and readability, we include in Section 4 a self-contained exposition of our construction.

We would like to thank Rachel Roberts for the many helpful conversations and useful advice and suggestions she gave throughout the writing of this paper. We would also like to thank Indira Chatterji for her interest in this class of partially ordered groups.

\section{Background information}

An {\it order tree} $T$ \cite{GO89} is a set $T$ together with a collection $\mathcal S$ of linearly ordered subsets called {\it segments}. If $\sigma$ is a segment then $-\sigma$ denotes the same subset with reverse order. The segments
satisfy :
\begin{itemize}
\item  Each segment $\sigma$ has distinct least and greatest elements, which we will denote by $i(\sigma)$ and $f(\sigma)$ respectively. (We also write $\sigma = [i(\sigma), f(\sigma)]$.)
\item If $\sigma$ is a segment, so is $-\sigma$.
\item A closed nondegenerate (i.e., containing more than one element) subinterval of a segment is a segment.
\item Any two elements of T can be joined by a sequence $\sigma_1,...,\sigma_k$ of segments such that $f(\sigma_j) = i(\sigma_{j+1})$ for all j.
\item Given a cyclic word $\sigma_0 \sigma_1 ...\sigma_{k-1}$ (where $i(\sigma_j) = f(\sigma_{j+1})$ for all $j, 0\leq j\leq k-2$, and cyclic means  $f(\sigma_{k-1}) = i(\sigma_0)$), there is a subdivision of the $\sigma$'s $\rho_0...\rho_{n-1}$ so that when adjacent pairs $(\rho)(-\rho)$ are cancelled, we have the trivial word.
\item If $f(\sigma_1) = i(\sigma_2) = \sigma_1 \cap \sigma_2$, then $\sigma_1 \cup \sigma_2$ is a segment.
\end{itemize}
An {\it $\mathbb{R}$-order tree} \cite{GK1} is an order tree satisfying also:
\begin{itemize}
\item Each segment is order isomorphic to a closed interval in $\mathbb{R}$.
\item $T$ is a countable union of segments.
\end{itemize}
\noindent
An {\it orientation} of an order tree is a choice of subset $\mathcal S_+ \subset \mathcal S$ such that
\begin{itemize}
\item $\mathcal S_+\cap (-\mathcal S_+) = \emptyset$, where $-\mathcal S_+ = \{ -\sigma|\sigma\in \mathcal S_+\}$.

\item A closed nondegenerate subinterval of a segment in $\mathcal S_+$ is in
$\mathcal S_+$.

\item Any two elements of T can be joined by a sequence $\sigma_1,...,\sigma_k$ of segments in $\mathcal S_+\cup (-\mathcal S_+)$ such that $f(\sigma_j) = i(\sigma_{j+1})$ for all j.

\item If $\sigma_1, \sigma_2\in\mathcal S_+$, and $f(\sigma_1)=i(\sigma_2)=\sigma_1\cap\sigma_2$, then $\sigma_1\cup\sigma_2\in\mathcal S_+.$
\end{itemize}
Since there are no nontrivial cyclic words, orientations always exist.

\begin{rmk}For simplicity of exposition, we will take all order trees to be $\mathbb R$-order trees.  Thus, unless otherwise noted, any order tree will be assumed to satisfy the above two $\mathbb R$ axioms.
\end{rmk}

Some special cases of order trees include $\mathbb{R}$-trees with countably many branch points and simply-connected (not necessarily Hausdorff) 1-manifolds. In the case of a simply-connected 1-manifold there are exactly two possible orientations, whereas for a general order tree there may be many more.  Any nontrivial orientation preserving group action on an oriented $\mathbb{R}$-order tree canonically induces an action on a related simply connected 1-manifold. The full details of the construction of this 1-manifold (a Denjoy blow-up of the original) appear in section 5 of \cite{RSS}, but for ease in reference later on we summarize the construction here, including some details.

First, we recall the notion of {\it incidence} for order trees.  We remark that the definition of incidence or branching degree extends to arbitrary order trees.  Fix an orientation on $T$ and let $x\in T$. Define an equivalence relation $\approx_f$ on the set $S(x,f)=\{\sigma\in \mathcal S_+|f(\sigma)=x\}$ by $\sigma_1\approx_f\sigma_2$ if and only if both $f(\sigma_1)=f(\sigma_2)=x$ and  $|\sigma_1\cap\sigma_2|>1$. For each $\sigma\in S(x,f)$, let $r_\sigma = \{\tau \in S(x,f)|\tau\approx_f\sigma\}$ and call $r_\sigma$ an {\it incoming ray} at $x$. Let $R(x,f)=\{r_\sigma|\sigma\in S(x,f)\}$. Call $n_f(x)=|R(x,f)|$ the {\it in degree} at $x$. Similarly, define an equivalence relation $\approx_o$ on the set $S(x,o)=\{\sigma\in \mathcal S_+|i(\sigma)=x\}$ by $\sigma_1\approx_o\sigma_2$ if and only if both $i(\sigma_1)=i(\sigma_2)=x$ and $|\sigma_1\cap\sigma_2|>1$. For each $\sigma\in S(x,o)$, let $r_\sigma = \{\tau \in S(x,o)|\tau \approx_o\sigma\}$ and call $r_\sigma$ an {\it outgoing ray} at $x$. Let $R(x,o)=\{r_\sigma|\sigma\in S(x,o)\}$. Call $n_o(x)=|R(x,o)|$ the {\it out degree} at $x$. We say that a segment $\sigma$ is incident to $x$ if $\sigma\in S(x,o)\cup S(x,f)$, and we say that a ray $r_\sigma$ is incident to $x$ if $r_\sigma\in R(x,o)\cup R(x,f)$.  Call $x\in T$ {\it regular} if $n_o(x)=n_f(x)=1$. Call $x\in T$ a {\it branch point} if it is not regular, and let $\mathcal B$ denote the set of branch points of $T$. Note that if $\mathcal B=\emptyset$, then $T$ can also be given the structure of a simply connected 1-manifold.

Now consider any $x\in \mathcal B$. If the out-degree $n_o(x)=0$ (in-degree $n_f(x) = 0$), call $x$ a {\it sink} (respectively, {\it source}).  If $n_o(x)=1$ and $n_f(x)>1$, call the single element $r_\sigma\in R(x,o)$ the {\it distinguished ray} at $x$. Symmetrically, if $n_f(x)=1$ and $n_o(x)>1$, call the single element $r_\sigma\in R(x,f)$ the {\it distinguished ray} at $x$.

\begin{lemma}[see Lemma 5.9 of \cite{RSS}]\label{l:distray}

Let $T_0$ be an oriented order tree such that at every $x\in\mathcal B$, there is a distinguished ray. Then any nontrivial orientation-preserving action on $T_0$ canonically induces a nontrivial orientation preserving action on a related oriented simply connected 1-manifold $T'$.
\end{lemma}

\begin{proof}
The 1-manifold $T'$ is obtained from the order tree $T$ by blowing up each branch point of $T$ into a set of endpoints for each ray except the distinguished ray, which is left open. See \cite{RSS} for full details.
\end{proof}

\begin{prop}[see Proposition 5.10 of \cite{RSS}]\label{P:tprime}
Any nontrivial ori-entation-preserving action on an oriented order tree $T$ canonically induces a nontrivial orientation preserving action on a related oriented simply connected 1-manifold $T'$.
\end{prop}

\begin{proof}

Again, a full proof appears in \cite{RSS}. At each branch point with in-degree and out-degree greater than 1, a linear Denjoy blowup (as in Definition 9.1 of \cite{RSS}) is performed to create a tree $T_0$ in which all branch points are sinks, sources, or have distinguished rays. A distinguished ray is then added at each sink and source, and then Lemma~\ref{l:distray} is applied to obtain the 1-manifold $T'$.
We show that any nontrivial orientation-preserving action on an oriented $\mathbb{R}$-order tree $T$ canonically induces a nontrivial orientation preserving action on an oriented $\mathbb R$-order tree $T_0$ such that at every $x\in\mathcal B$, there is a distinguished ray. Lemma \ref{l:distray} then applies.
\end{proof}

In the construction described above, there is a natural map $\varphi:T' \rightarrow T$.  This map collapses to a point each segment added during the Denjoy blow-up of $T$ to $T_0$, collapses to the sink or source each segment added to these points to give them distinguished rays and collapses all the points in each set $\{x_{r_\sigma}\}$ added for each distinguished ray $\hat{r}_x$.  We will need to keep track of an important implicit subtree of $T'$ which maps surjectively onto $T$, called the core of $T'$.

\begin{definition}
$\hat T$, the \emph{core} of $T'$, is the subset of $T'$ consisting of all points not in the union of open rays which are added to $T$ at the last step in forming $T_0$ (in order to transform branch points which are sources or sinks into branch points with a distinguished ray).

\end{definition}

Next we recall from \cite{RS2} some of the basics about the structure of order trees.

\begin{definition} Let $T$ be an order tree. A \emph{path} from $x \in T$ to $y \in T$ is a sequence of segments $\sigma_1 \cdots \sigma_n$ with $f(\sigma_i)=i(\sigma_{i+1})$ for $1 \le i < n$ and $i(\sigma_1)=x$ and $f(\sigma_n)=y$.
\end{definition}

\begin{definition} A \emph{standard geodesic} from x to y is a path $\sigma_1 \cdots \sigma_n$ from x to y satisfying:

\begin{itemize}

\item $\sigma_i \cap \sigma_j= \emptyset$ if $|i-j| > 1$

\item $\forall j$, either $\sigma_j \cap \sigma_{j+1} = i(\sigma_{j+1})=f(\sigma_j)$ or $\sigma_j \cap \sigma_{j+1} =(i(\sigma_j),f(\sigma_j)]= (f(\sigma_{j+1}),i(\sigma_{j+1})]$
\end{itemize}
\end{definition}

Let $S=\{(\sigma, \tau)| (i(\sigma),f(\sigma)]=(f(\tau),i(\tau)]$ and $i(\sigma)\neq f(\tau)\} $, where $\sigma$ and $\tau$ are segments. We define a relation on $S$ as follows: let $([x,z],[z,y]) \equiv ([x,z'], [z',y])$ if $\exists r \in (x,z] \cap (x,z']$ so that $[x,z]=[x,r] \cup [r,z], [z,y]=[z,r]\cup[r,y], [x,z']=[x,r]\cup[r,z'],[z',y]=[z',r] \cup[r,y]$ (where the segments $[r,z]$ and/or $[r,z']$ are understood to be empty if $r=z$ or $r=z'$).  Then $\equiv$ is an equivalence relation.

\begin{definition}
A {\it cusp} is an equivalence class of pairs of segments in $S$ under the above equivalence relation $\equiv$.
\end{definition}

\begin{notation}
Note that a cusp represented by a pair $([x,z],[z,y]) \in S$ is determined by the pair of points $x$ and $y$, for by axiom 5 in the definition of order tree, any other pair in $S$ of the form $([x,z'],[z',y])$ must be in the same equivalence class. Hence we will use the symbol $[x,y]^c$ to denote the cusp.
In this situation, we refer to the points $x$ and $y$ as {\it cusp points}.

\end{notation}

Then we have the following existence theorem.

\begin{thm}[see~\cite{RS2}, Theorem 3.4]
Let $T$ be an order tree. Given $x$ and $y$ $\in T$, $\exists$  a standard geodesic from x to y.
\end{thm}

\begin{rmk}
In particular, if a path satisfies the first condition in the definition of a standard geodesic then no three segments in the path have a nonempty intersection.  Hence each segment in a standard geodesic is a part of a representative of at most one cusp, or in other words it is never the case that $(\sigma_j,\sigma_{j+1})$ and $(\sigma_{j+1},\sigma_{j+2})$ are both representatives of cusps.
\end{rmk}

Standard geodesic paths are not unique, but the lack of uniqueness as a set of points is all due to the lack of uniqueness in the representation of cusps as two segments. We make this more precise with the following definition.

\begin{definition}
Let $\gamma$ be a standard geodesic from $x$ to $y$. Define $GS_{(x,y)}$, the {\it geodesic spine} of $\gamma$,  to be the union of all segments of $\gamma$ which do not, together with an adjacent segment, give a representative of a cusp, together with all of the cusp points of $\gamma$.
\end{definition}

Although $GS_{(x,y)}$ is defined as a set, it has a natural linear order inherited from the geodesic $\gamma$. We will sometimes abuse language and call the geodesic spine a path, even though it has gaps between cusp points. When a pair of segments representing a cusp $[p_1,p_2]^c$ is on $\gamma$, it is only $p_1$ and $p_2$ which are on $GS_{(x,y)}$. However, we will again sometimes abuse language and say that $[p_1,p_2]^c$ is on $GS_{(x,y)}$ to stress the fact that $p_1$ and $p_2$ are cusp points.

As the notation suggests, the geodesic spine of $\gamma$ depends only on the endpoints of the geodesic. It is independent of the particular choice of standard geodesic $\gamma$. The uniqueness of the geodesic spine of a standard geodesic between $x$ and $y$ will follow from the following theorem.

\begin{thm}[see~\cite{RS2}, Theorem 3.6]
The geodesic spine of a standard geodesic from $x$ to $y$ is the intersection of all paths from $x$ to $y$.
\end{thm}

The standard notions of trivial actions (ones with global fixed points) and minimal actions (those with no proper invariant subtree) must be modified slightly to suit the order tree situation.

\begin{definition}\label{D:genfixedpoint}
If $g \in G$ where the group $G$ acts on an order tree $T$, we say that $g$ has a generalized fixed point $x\in T$ if  $gx$ is not separated from $x$ (every segment containing $x$ intersects every segment containing $gx$).
\end{definition}

An action of $G$ on $T$ is nontrivial if there is no point $x\in T$ which is a generalized fixed point for every element of $G$.  Another type of action which is trivial in spirit is an action with a unique fixed end. In the case of an $\mathbb{R}$-tree, an end can be defined as an equivalence class of rays, where a ray is an embedding of $[0, \infty)$.  We adapt the definitions slightly to the case of a (possibly not Hausdorff) order tree.

\begin{definition}\label{D:ray2}
A {\em ray} in an order tree $T$ is a subset $\rho$ of $T$ that can be written as an infinite increasing union of geodesic spines,
$$ \rho = \bigcup_{i=1}^\infty GS_{(x,x_i)},$$
where $GS_{(x,x_i)} \subset GS_{(x,x_j)}$ if $i<j$.  The ray $\rho$ is said to be {\em infinite} if the sequence $\{x_i\}$ does not converge in $T$.
\end{definition}

\begin{definition}\label{D:endpoint}
If the ray $\rho$ is not infinite, we call a limit of the sequence $\{x_i\}$ an {\em endpoint} of $\rho$.  If $T$ is not Hausdorff, a finite ray may have multiple endpoints, but all endpoints of the same ray are non-separable from each other.
\end{definition}

\begin{lemma}\label{L:ray.is.gs}
If $\rho$ is a finite ray starting at the point $y$ of an order tree $T$ and $a$ is an endpoint of $\rho$ then $\rho = \gs[a]y - \{a\}$.
\end{lemma}

\begin{proof}
Let $\rho = \cup_{i=1}^\infty \gs[y]{y_i}$.  First suppose that $\alpha \in \gs[a]y$ and $\alpha \neq a$.  If $\alpha = y$ then $\alpha = y \in \rho$.  If $\alpha \neq y$ then $\alpha$ separates $y$ and $a$. All but finitely many $y_i$ lie in the neighborhood of $a$ consisting of the connected component of $T- \{\alpha\}$ that contains $a$. So there is a $j$ such that $\alpha$ does not separate $y_j$ from $a$, which implies that   $\alpha \notin \gs[a]{y_j}$. Hence  $\alpha \in \gs[y]{y_j}$, which implies that $\alpha \in \rho$.  Thus $\gs[a]y -\{a\} \subset \rho$.

For the other inclusion, suppose that $y_i \notin \gs[a]y$.  Then $a$ and $y$ lie in the same component of $T - \{y_i\}$.  But for $j > i$, $y_i$ separates $y$ from $y_j$, so $y_j$ does not lie in the connected component of $T - \{y_i\}$ that contains $a$.  Since this is true of all $y_j$ with $j > i$, the sequence $\{y_i\}$ cannot have $a$ as a limit point contradicting the fact that $a$ is an endpoint of $\rho$.  So for each $i$, $y_i \in \gs[a]y$ and $\rho \subset \gs[a]y$.  Now, $a\notin \rho$ for otherwise $a$ would have to lie in $\gs[y]{y_i}$ for some $i$.  In this case either $a=y_i$, which is nonsense, or only finitely many $y_j$ would lie in the connected component of $T-\{y_i\}$ containing $a$, namely $y_j$ with $j<i$.  Therefore, $\rho \subset \gs[a]y -\{a\}$.
\end{proof}

There is an equivalence relation on infinite rays in $T$ given by $\rho_1 \approx \rho_2$ if the intersection $\rho_1 \cap \rho_2$ contains an infinite ray.  Note that since geodesic spines are unique, whenever $\rho_1 \cap \rho_2$ contains a infinite ray, the two rays $\rho_1$ and $\rho_2$ eventually coincide.  Additionally, if $\rho_1 \not\approx \rho_2$, then $\rho_1 \cap \rho_2$ is contained in a finite geodesic spine, and the two rays $\rho_1$ and  $\rho_2$ eventually separate.

\begin{definition}\label{D:end}
An {\em end} of the order tree $T$ is an equivalence class of infinite rays in $T$.
\end{definition}

Now, $G$ acts on the set of infinite rays and for any $g\in G$, $\rho_1 \approx \rho_2$ if and only if $g\rho_1 \approx g\rho_2$.  Therefore, $G$ acts on the set of ends of $T$.

And finally, the standard notions of invariant subtrees and minimal actions must be adapted to this situation.

\begin{definition}
If $T$ contains a $G$-invariant subset $T^{\prime}$, with the property that for any two points $x,y \in T'$, $GS_{(x,y)}\subseteq T'$ , we call $T^{\prime}$ an invariant {\it implicit subtree}. Of course, an invariant subtree is one special case of an invariant implicit subtree. We call $T^{\prime}$ an invariant {\it implicit line} if $T^{\prime}$ admits a total ordering (without a greatest or a least element) with the following property: Choose any $x,y \in T^\prime$ and let $[x,y]$ denote the interval of $T^\prime$ determined by the total ordering. Then $GS_{(x,y)}=[x,y]$ and furthermore, the natural ordering on $GS_{(x,y)}$ agrees with the total order on $[x,y]$.
\end{definition}

\begin{definition}Let a group $G$ act on an order tree $T$. The action is {\it minimal} if $T$ contains no proper invariant implicit subtree.
\end{definition}

Analogous to the case for Hausdorff trees we have,

\begin{lemma}\label{L:twofixedends}
If $G$ acts on the order tree $T$ fixing two distinct ends, then $G$ fixes an implicit line $l$, the ends of which are fixed by $G$.
\end{lemma}

\begin{proof}
Suppose that $G$ fixes the ends $\epsilon_1$ and $\epsilon_2$.  Choose rays $\rho_1 = \cup_{i=1}^\infty \gs[a]{a_i}$ and $\rho_2 = \cup_{i=1}^\infty \gs[b]{b_i}$ representing $\epsilon_1$ and $\epsilon_2$ respectively.  Now, $\rho_1$ and $\rho_2$ may coincide for a time, but since $\epsilon_2 \neq \epsilon_2$ the rays $\rho_1$ and $\rho_2$ eventually separate.  Choose points $a_m$ and $b_n$ after the point of separation.  Then
\begin{equation*}
l := \left(\bigcup_{i=m}^\infty \gs[a_m]{a_i}\right) \cup \left( \gs[a_m]{b_n} \right) \cup \left( \bigcup_{i=n}^\infty \gs[b_n]{b_i} \right)
\end{equation*}
is an implicit line.  Note that the ends of $l$ are fixed by $G$ because one represents $\epsilon_1$ and the other represents $\epsilon_2$.

To see that $l$ is fixed by $G$, let $p\in l$ and $g\in G$.  Then $p$ separates $l$ into two infinite rays,
\begin{eqnarray*}
\alpha_1 & = & \bigcup_{i=1}^\infty \gs[p]{p_i} \\
\alpha_2 & = & \bigcup_{i=1}^\infty \gs[p]{q_i}
\end{eqnarray*}
with $\alpha_i$ representing $\epsilon_i$.  Since $G$ fixes $\epsilon_1$ and $\epsilon_2$, $g\alpha_1$ eventually coincides with  $\alpha_1$ and $g\alpha_2$ eventually coincides with $\alpha_2$.  Choose $p_m \in g\alpha_1 \cap \alpha_1$ and $q_n \in g\alpha_2\cap \alpha_2$.  Let $\sigma_1$ be the initial segment $\gs[p]{p_m}$ and $\sigma_2$ be the initial segment $\gs[p]{q_n}$.  Then $\sigma_1\cap\sigma_2 = \{p\}$, so that $g\sigma_1 \cap g\sigma_2 = \{gp\}$.  Since $g\sigma_i$ is an initial segment of $g\alpha_i$, the initial segments $\gs[gp]{p_m}$ and $\gs[gp]{q_n}$ intersect only in the one point set $\{gp\}$.  Therefore, $\gs[p_m]{gp} \cup \gs[gp]{q_n}$ is the geodesic spine $\gs[p_m]{q_n}$.  But $\gs[p_m]{q_n}$ is contained in $l$.  So
$$gp \in \gs[p_m]{gp} \cup \gs[gp]{q_n} = \gs[p_m]{q_n} \subset l,$$
showing that $gp \in l$ for any $p \in l$ and $g \in G$.  Therefore, $G$ fixes the implicit line $l$.
\end{proof}

\section{Partial orders}\label{s:order}

The points of an oriented simply connected 1-manifold $T$ are often considered to be partially ordered in a natural way, by declaring two points to be comparable if and only if they lie in submanifold which is homeomorphic to $\mathbb{R}$ by order preserving homeomorphism (the smaller of a comparable pair of points is determined by the orientation on the manifold). There is an extension to this partial order, which is spelled out in definition 4.4 of \cite{RSS}, which is also naturally inspired by the orientation on the manifold.  Moreover, it is the maximal extension which is guaranteed to be preserved by all orientation preserving homeomorphisms of the manifold. We review this order here.  If $x \in T$, let $I_x$ be an open set in $T$ containing $x$ which is homeomorphic (as an oriented manifold) to $\mathbb{R}$. Then there is a total order on the points of $I_x$ which is induced by the homeomorphism. Let $I_x^+$ be the set of elements of $I_x - \{x\}$ which are greater than $x$, and let $I_x^-=I_x-(I_x^+ \cup \{x\})$. Then $T-\{x\}$ has exactly two connected components, and we let $x^+$ be the component containing $I_x^+$ and let $x^-$ be the component containing $I_x^-$. Then a partial order on the points of $T$ is given by $x \leq y \iff y^+ \subseteq x^+$.

\begin{notation} For an arbitrary poset $S$, and for $x,y \in S$, if $x$ and $y$ are incomparable, write $x \sim y$. If so, and they have a common upper bound, write $x \sim_u y$, and if they have a common lower bound write $x \sim_l y$. In a general partial order, a pair of elements could have either no common bounds, just a common upper bound, just a common lower bound, or both types of common bounds. \end{notation}

We observe that the poset of the points on a simply connected 1-manifold satisfies the following properties, which a general poset may or may not satisfy:

\begin{definition}\label{p:bounds} A poset $S$ is \emph{strongly connected} if for any pair of incomparable elements of the poset $S$, say $x \sim y$, either $x \sim_u y$ or $x \sim_l y$.
\end{definition}

We remark that this is, as the name suggests, stronger than the standard definition of a connected poset. A poset is generally considered to be connected if given any pair of elements $x,y\in S$, there is a finite sequence $x=a_0, a_1, \ldots , a_n=y$ in $S$ such that $a_i$ is comparable to $a_{i+1}$ for all $0 \leq i \ < n$ (see \cite{KM}). Using this language, a totally ordered set is a connected poset in which a sequence can always be chosen to have $n=1$, and a strongly connected poset is a connected poset where a sequence can always be chosen with $n\leq 2$.

\begin{definition}\label{p:nobadtriads}
A poset $S$ is \emph {acyclic} if $ \forall x,y,z \in S$, $x \sim_u y$ and $x \sim_l z \Rightarrow z>y$.
\end{definition}

\begin{definition}\label{simplyconnected}
A partially ordered set will be called \emph{simply connected} if it is both strongly connected and acyclic.
\end{definition}

Note that in particular, the partial order on the points of an oriented simply connected 1-manifold is a simply connected partial order.

\begin{rmk}
If two elements in a simply connected partially ordered set $S$ are not comparable, then acyclicity implies that it is not possible for the pair to have both common upper and lower bounds, and strong connectivity implies that the pair must have one type of common bound. Hence, for $x \neq y$ elements of $S$, exactly one of the four possible relationships hold, namely one of $x<y$, $x>y$, $x \sim_l y$, and $x \sim_u y$ holds.
\end{rmk}

Our goal is to characterize the partial orders which groups acting (nontrivially and without fixing a unique end) on oriented order trees must admit, and it turns out that simply connected partial orders are too restrictive.  For a complete characterization, the notion of a formal extension is necessary.

\begin{definition}\label{D:extension}
Let $S$ be an acyclic partially ordered set. A \emph{formal extension} of the partial order is the set $S$ along with the partial order, with an additional structure as follows. Each pair $x \sim y$ which has neither a common upper or lower bound is formally assigned exactly one of the two types $\sim_u$ or $ \sim_l$, so that for each pair $x \neq y$, exactly one of the four possible formal relationships holds ($x<y, x>y, x \sim_l y, x\sim_u y$). The formal extension is said to be a \emph{simply connected extension} if the resulting set of formal relationships, which by construction is strongly connected, also formally satisfies acyclicity. \end{definition}

\begin{figure}
\begin{center}
\psfrag{T}{$T$}
\psfrag{T'}{$T'$}
    \includegraphics[width = 3in]{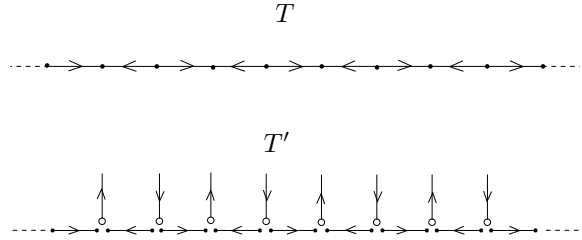}
    \caption{An oriented order tree structure on $\mathbb R$ and associated $1$-manifold}\label{F:dihedral}
\end{center}
\end{figure}

The necessity for this definition is seen in the infinite dihedral group.  This group does not admit a nontrivial partial order which is both strongly connected and acyclic. However, it does have a partial order with a simply connected formal extension. This order comes naturally from an action of the group on the real line, where the real line is viewed as an order tree $T$ with maximal positively oriented segments $\sigma_i =[i,i+1]$, $i \in \mathbb{Z}$  with $o(\sigma_i)=i$ and $f(\sigma_i)=i+1$ if $i$ is even, and $o(\sigma_i)=i+1$ and $f(\sigma_i)=i$ if $i$ is odd, as shown in Figure~\ref{F:dihedral}.  Positive translations shift $T$ two units to the right, and reflections reflect about the integer points.  Note that the action of the dihedral group on $T$ is minimal, while the action on the associated non-Hausdorff $1$-manifold $T'$ stabilizes any submanifold consisting of the horizontal geodesic spine together with all the translates of any open proper subinterval at the bottom of a single vertical edge. Therefore, the action of the dihedral group on $T'$ is not minimal, nor does it contain a minimal $G$-invariant submanifold. The group is embedded naturally in $T'$ by identifying it with the orbit of a point with trivial stabilizer, and the simply connected partial order on $T'$ induces a partial  order on this orbit, and hence the group, which is acyclic but is not even connected, let alone strongly connected. However, the embedding in $T'$  gives a recipe for a formal extension of the order structure which is simply connected.


Note that the empty partial order on any set, by which we mean the order in which nothing is comparable to anything else, always has a trivial simply connected extension by just declaring all pairs of elements to be the same incomparability type. Hence we define:
\begin{definition}\label{d:simplyconnectedextensiontrivial}
A simply connected extension of a partial order is \emph{trivial} if at least one pair of elements is non-comparable, and all such non-comparable pairs are of a single type.
\end{definition}

\begin{rmk}
A partial order which is already simply connected has a simply connected extension -- the one in which no additional relationships are assigned. Such a partial order may or may not be trivial in the above sense. In addition, note that there are certainly groups which cannot admit nontrivial extensions; for example, a torsion group. More interesting examples from the point of view of actions on trees appear in Section 5.
\end{rmk}

We record next an elementary lemma about these partial orders.

\begin{lemma}\label{L:trivial}  If $x,y,z \in S$, where $S$ is a set with a simply connected extension of a partial order, then:
\begin{itemize}
\item $x \sim_l y$ and  $y<z \Rightarrow x \sim_l z$
\item $x \sim_u y$ and  $y>z \Rightarrow x \sim_u z$.
\end{itemize}
\end{lemma}

\begin{proof} First suppose that $x \sim_l y$ and  $y<z$. So $x$ and $y$ are incomparable, and have no common upper bound. Recall that $x$ and $z$ always satisfy exactly one of the four relationships. We show that $x \sim_l z$ by eliminating the other three possibilities. First, $x<z$ implies that $z$ is a common upper bound for $x$ and $y$, which is impossible. But $x>z \Rightarrow x>y$, which is not true. Now if $x \sim_u z$, then by acyclicity it must be true that $y>z$, which is once again a contradiction. Hence it must be that $x\sim_l z$. The proof of the second statement is similar.
\end{proof}

We will be interested in the special case where the partially ordered set is a group and the order is left-invariant. In that setting, if $g\sim h$ in the original partial order, and the pair have neither an upper nor a lower common bound, then of course no pair of translates $fg \sim fh$ will either. We will, of course, only be interested in extensions which are left-invariant.

\begin{definition}
Let $G$ be a group with a left-invariant acyclic partial order. We say it has a left-invariant simply connected extension if it has a simply connected extension as a poset which is left-invariant in the sense that $g\sim_u h \Rightarrow fg\sim_ufh~~\forall f \in G$, and similarly for $\sim_l$ .
\end{definition}

For use in the later sections, we investigate the properties of an acyclic partial order on a set $S$  which we assume to have a simply connected extension. In what follows, we are referring always to the full set of comparability relationships in a chosen simply connected extension of the given partial order. Everything works, of course, in the special case where the partial order is already simply connected and doesn't need to be extended.

First we define the notion of betweenness inspired by the geodesic spines in $\mathbb{R}$-order trees.

\vspace{.1 in}
\begin{definition}{(Betweenness):}\label{D:betweenness}
Definition: Given two points $a$ and $c$ with $c \neq a$, we say that $b$ is between $a$ and $c$ if any of the following three conditions holds:
\begin{enumerate}
\item  In the case that $a<c$, either $a<b<c$ or $a \sim_u b$ and $b \sim_l c$.
\item  In the case that $a \sim_l c$, either $a \sim_l b$ and $b<c$ or $c \sim_l b$ and $b<a$.
\item  In the case that $a \sim_u c$, either $a \sim_u b$ and $b>c$ or $c \sim_u b$ and $b>a$.
\end{enumerate}
\end{definition}

\begin{notation} For any $a \neq b$, $a,b \in G$, let $$B_{a,b}=\{a,b\} \cup \{x\mid x~~is~~between~~a~~and~~b\}.$$
\end{notation}

The following theorem, which can be easily proven by careful case-by-case analysis, records some basic properties of between sets.

\begin{thm} Let $a,b,c,d \in S$, where $S$ is a partially ordered set with a simply connected extension. Then:
\begin{enumerate}
\item $\forall a,b,c, B_{a,b} \subset B_{a,c} \cup B_{c,b}$.

\item $c \in B_{a,b} \iff B_{a,b} = B_{a,c} \cup B_{c,b}$.

\item If $c \in B_{a,b}$, then $B_{a,c} \cap B_{c,b} = \{c\}$.

\item If $b \in B_{a,c}$ and $c \in B_{b,d}$, then $b,c \in B_{a,d}$.
\end{enumerate}
\end{thm}

From this theorem, it is easy to see that these sets $B_{a,b}$ come with a natural total order on them. The only arbitrariness is in declaring whether $a$ is the least or greatest element.

\begin{definition}
 If $x$ and $y \in B_{a,b}$, $x\preceq y$ if $x \in B_{a,y}$ and $x\succeq y$ if $x \in B_{y,b}$. (Of course, here $a$ is being considered least element and $b$ the greatest).
\end{definition}

\begin{cor}
The order on $B_{a,b}$ is a total order.
\end{cor}

\begin{definition} A path from $x$ to $y$ is a finite union of $n$ sets $B_{s_i,t_i}$ with $t_i=s_{i+i}$, $s_1=x$, $t_n=y$.

\end{definition}

\begin{cor} $B_{x,y}$ is the intersection of all paths from $x$ to $y$.

\end{cor}

Note that if the set $S$ is the set of points on an oriented simply connected 1-manifold with the order described earlier, then $B_{a,b}=GS_{(a,b)}$. Since these geodesic spines can always be expressed as finite unions of oriented segments in order trees, we suspect that our extensions of partial orders will require a finiteness condition on the paths $B_{a,b}$. To articulate this, we define a relation on each set $B_{a,b}$.

\begin{definition}
For $x,y \in B_{a,b}$,  $xO_{a,b}y$ if and only if $B_{x,y}$ is totally ordered in the original order on the set $S$.
\end{definition}

\begin{lemma}
The relation  $O_{a,b}$ is an equivalence relation on $B_{a,b}$.
\end{lemma}

\begin{proof} $O_{a,b}$ is clearly reflexive and symmetric. Transitivity follows from the properties of the extension of the order on $S$.
\end{proof}
It is now clear that the equivalence classes are themselves totally ordered in the original order, and that if $x$ and $y$ are in the same equivalence class, then all points between $x$ and $y$ are also in that equivalence class. 

\begin{definition}\label{D:finiteswitches}
 We say that the simply connected extension of a partial order on a set $S$ is \emph{rectifiable} if $\forall a,b \in S, B_{a,b}$ is a finite union of equivalence classes under the relation $O_{a,b}$.
\end{definition}

Partial orders with rectifiable simply connected extensions are naturally inspired by the partial orders on simply connected 1-manifolds, and they satisfy algebraic properties commonly found in the theory of partially ordered groups. We now discuss a few of these. Firstly, it is obvious that if a group has such a partial order, then so does any subgroup (though the order restricted to the subgroup may be trivial). Furthermore, these orders also pass to appropriate quotient groups. It is a standard result in the theory of partially ordered groups that if $G$ has a left invariant partial order, and $H$ is a normal, convex subgroup, then a partial order is naturally induced on the quotient group $G/H$ (see section1.6.3 of \cite{KM}). In the case of a simply connected partial order, or more generally a simply connected extension, the order on the quotient group is of the same special type as long as we require a stronger version of convexity for the subgroup, obtained by replacing the notion of totally ordered sets by between sets.

\begin{definition}\label{D:completelyconvex}
Let $G$ be a group with a left-invariant partial order with left-invariant simply connected extension. Then a subgroup $H$ is called \emph{completely convex} if for any pair of elements $h_1,h_2 \in H$, $B_{h_1,h_2}\subset H$.
\end{definition}

Then we have:

\begin{thm}\label{T:quotientgetsorder}Let G be a group with a left-invariant partial order with left-invariant simply connected extension. Let $H$ be a completely convex normal subgroup. Then $G/H$ also admits a left-invariant partial order with left-invariant simply connected extension. Furthermore, if the extension is rectifiable, then the extension induced on $G/H$ is also rectifiable.
\end{thm}

\begin{proof}
First we define four possible relationships between cosets in $G/H$. Given cosets $g_1H \neq g_2H$, define

\begin{enumerate}
\item $g_1H<g_2H$ if $\exists h \in H$ such that $g_1 < g_2h$.
\item $g_1H>g_2H$ if $\exists h \in H$ such that $g_1 > g_2h$.
\item $g_1H \sim_u g_2H$ if $\exists h \in H$ such that $g_1 \sim_u
g_2h$.
\item $g_1H \sim_l g_2H$ if $\exists h \in H$ such that
$g_1 \sim_l g_2h$
\end{enumerate}
Since $H$ is normal and the partial order is left invariant, if a pair of cosets satisfies one of the above conditions for a particular choice of coset representatives, then they will satisfy that condition with any choice of coset representatives.  Also, for any pair of distinct  cosets $g_1H$ and $g_2H$, either $g_1<g_2$, $g_2<g_1$, $g_1 \sim_u g_2$, or $g_1 \sim_l g_2$, since we have a full extension of a partial order in $G$.  So each pair of distinct cosets will be related in at least one of the above ways. In addition, it is clear that  these relationships are invariant under left multiplication by elements of $G$.

We now prove that the relationships defined above satisfy the following property:

\begin{property}\label{prop2like} For $g_1, g_2, g_3 \in G$,

\begin{enumerate}
\item $g_1H < g_2H$ and $g_2H<g_3H$ imply that $g_1H < g_3H$.
\item $g_1H \sim_u g_2H$ and $g_2H \sim_l g_3H$ imply that $g_3H <
g_2H$.
\item $g_1H \sim_u g_2H$ and $g_3H < g_2H$ imply that $g_1H
\sim_u g_3H$
\item $g_1H \sim_l g_2H$ and $g_2H < g_3H$ imply that $g_1H \sim_l
g_3H$
\end{enumerate}
\end{property}
\begin {proof}
The proofs of the four statements are virtually identical. In each case, the normality of $H$ gives the desired relationship as long as $g_1H \neq g_3 H$, and then the complete convexity of $H$ implies that the two cosets cannot in fact be the same.  We provide the explicit argument for the first statement. If $g_1H < g_2H$ and $g_2H<g_3H$, then $g_1 < g_2h$ and $g_2 < g_3h'$ for some $h,h' \in H$. So by the normality of $H$, $g_2h < g_3h''$ for some $h'' \in H$. But then $g_1 < g_3h''$, which implies that $g_1H < g_3H$, the desired conclusion, as long as $g_1H \neq g_3H$. But $g_3H=g_1H$ implies that $g_1^{-1}g_3 \in H$. Notice that translating the relationships above by $g_3^{-1}$, we see that $g_3^{-1}g_1 < g_3^{-1}g_2h < h''$. But this implies that $g_3^{-1}g_2h \in B_{g_3^{-1}g_1,h''}$, which in turn is contained in $H$ by the complete convexity of $H$. Hence $g_3^{-1}g_2 \in H$, or equivalently $g_3H=g_2H$, which is ruled out by the assumption of $g_2H < g_3H$.
\end{proof}

It is a formal consequence of Property \ref{prop2like}  that a given pair of distinct cosets can satisfy at most one of the four possible relationships. To prove this, there are four possibilities which must be eliminated for the pair $g_1H \neq g_2H$.
\begin{enumerate}
\item  Suppose $g_1H < g_2H$ and $g_2H < g_1H$. Then by the first part of Property \ref{prop2like}, $g_1H < g_1H$, which is impossible. \item  Suppose that $g_1H \sim_u g_2H$ and $g_2H \sim_l g_1H$. Then by the second part of Property \ref{prop2like}, $g_1H < g_1H$.

\item Suppose $g_1H \sim_u g_2H$ and $g_1H < g_2H$. Then by the third part of Property \ref{prop2like}, $g_1H \sim_u g_1H$, which is impossible.

\item Suppose $g_1H \sim_l g_2H$ and $g_2H < g_1H$. Then by the fourth part of Property \ref{prop2like}, $g_1H \sim_l g_1H$.
\end{enumerate}

At this point, we see that we have defined a partial order and a formal set of relationships $\sim_u$ and $\sim_l$ which satisfy all the properties of a simply connected extension. Finally, to ensure that this structure is really an extension of the partial order, we must check that the formal relationships $\sim_u$ and $\sim_l$ actually agree with any existing relationships based on common bounds. Namely, suppose that $g_1H$ and $g_2H$ are two distinct cosets which are not comparable according to the above definition of comparability, but do have a common upper (respectively lower) bound. We must show that in the above formal assignments, it is indeed the case that $g_1H \sim_u g_2H$ (respectively $g_1H \sim_l g_2H$). We argue the case of a upper bound. Suppose that $g_3H$ is this common upper bound. Then $g_1 \leq g_3h$ and $g_2\leq g_3h'$ for some $h,h' \in H$. Hence by normality of $H$, $g_2h'' \leq g_3h$ for some $h'' \in H$, so $g_1$ and $g_2h''$ share an upper bound in $G$. Now if $g_1$ and $g_2h''$ were comparable in $G$, then $g_1H$ and $g_2H$ would be comparable in $G/H$. Hence $g_1$ and $g_2h''$ are not comparable, and $g_1 \sim_u g_2h''$ in $G$. But this implies that $g_1H \sim_u g_2H$ in $G/H$. The argument for lower bounds is the same.

To complete the proof of Theorem \ref{T:quotientgetsorder}, we must show that if the original extension was rectifiable, then so is the induced extension of the order on $G/H$. We prove this by contradiction. Suppose that the extended order on $G/H$ is not rectifiable. Since the order is left-invariant, we may assume that for some $g \in G$, there are infinitely many equivalence classes in $B_{gH,H}$ under the relation $O_{gH,H}$. Then it follows that we may choose a sequence of cosets $\{g_0'H=gH,g_1'H,g_2'H, g_3'H, \ldots \}$ such that $g_i'H \neq H$, $g_i'H \in B_{g_{i-1}'H,H}$, $g_i'H \in B_{g_{i-1}'H,g_{i+1}'H}$, and $g_i'H$ not comparable to $g_{i-1}'H$ for all $i \geq 1$. From this we will inductively construct a sequence $\{g=g_0, g_1,g_2,\ldots \}$ and $h \in H$ satisfying $g_iH=g_i'H$, $g_i \in B_{g_{i-1},h} $, $g_i \in B_{g_{i-1},g_{i+1}}$, and $g_{i-1}$ not comparable to $g_i$ for any $i \geq 1$.  But this is impossible, for $\{g_i | i \geq 1 \}$ all lie in $B_{g,h}$, yet they all are in distinct equivalence classes under $O_{g,h}$, violating the fact that the extension of the order in $G$ is rectifiable. Before constructing the sequence, we state an easily verified lemma:
\begin{lemma}\label{l:btwnquotient}
Suppose that $gH \in B_{fH,kH}$ where $g,f,k \in G$. Then $gh' \in
B_{f,kh}$ for some $h,h' \in H$.
\end{lemma}

To construct the desired sequence, first note that $g_1'H \in B_{gH,H}$, so by Lemma \ref{l:btwnquotient} $g_1'h_1 \in B_{g,h}$ for some $h,h_1 \in H$. Set $g_1=g_1'h_1$, and then $g_1 \in B_{g_0,h}$ and $g_1$ is not comparable to $g_0$ as desired.  For the inductive step, suppose we have already constructed $\{g_1,g_2, \ldots, g_n \}$ with the desired properties. Since $g_{n+1}'H \in B_{g_nH,H}$, Lemma \ref{l:btwnquotient} yields $h_{n+1},h' \in H$ so that $g_{n+1}'h_{n+1} \in B_{g_n,h'}$. Again, let $g_{n+1}=g_{n+1}'h_{n+1}$. So we have already that $g_{n+1}H=g_{n+1}'H$, and also $g_{n+1} \in B_{g_n,h'}$.

Now by the complete convexity of $H$, $B_{h,h'} \subset H$, which implies that $g_{n+1} \notin B_{h,h'}$ since $g_{n+1}\notin H$. Hence it must be that $g_{n+1} \in B_{g_n,h}$ as desired, for if not, $g_{n+1} \notin B_{h,h'} \cup B_{g_n,h}$, which is nonsense given $g_{n+1} \in B_{g_n,h'} \subset B_{h,h'} \cup B_{g_n,h}$. In addition, it is clear that $g_n \in B_{g_{n+1},g_{n-1}}$ and $g_{n+1}$ and $g_n$ are incomparable. So we have inductively constructed the sequence needed to complete the proof.
\end{proof}

Next, in the theory of left-invariant partially ordered groups, the positive cone is an important tool, and there is an analogous structure for left-invariant simply connected extensions.  More precisely, in a left-invariant partially ordered group $G$, the {\em positive cone} is the subset $\mathcal P = \{g\in G \mid e<g\}$. The positive cone satisfies the two conditions,
\begin{eqnarray}
\PP \cap \PP^{-1} &=& \emptyset \label{E:ppure} \\
\PP \cdot \PP &\subset & \PP \label{E:ptimesp},
\end{eqnarray}
where $\PP^{-1}=\{x^{-1}|x\in \PP \}$.  The importance of the positive cone comes from the fact that any subset $\PP \subset G$ satisfying conditions~\ref{E:ppure} and~\ref{E:ptimesp} defines a left-invariant partial order on $G$ by the equation $g<h$ if and only if $g^{-1}h\in \PP$.  The partial order defined by $\PP$ is a total order if and only if $G = \PP \cup \PP^{-1} \cup\{e\}$.  Many other properties of partially ordered groups may be defined and studied via $\PP$.  For example a partially ordered group is {\em directed up} if any two elements of the group share an upper bound, which is equivalent to the condition,
\begin{eqnarray*}
G = \PP \cup \PP^{-1} \cup \left(\PP\cdot\PP^{-1}\right).
\end{eqnarray*}
Similarly, partially ordered groups with left-invariant simply connected extensions may be classified in terms of the positive elements and the elements sharing lower or upper bounds with the identity.  We have,

\begin{prop}\label{P:landu}
The group $G$ has a left-invariant partial order with simply connected extension if and only if there exist subsets $\PP, \LL$ and $\UU$ of $G$ such that,
\begin{enumerate}
    \item  $\PP \cap \PP^{-1} = \emptyset$, $\LL^{-1} = \LL$, $\UU^{-1} = \UU$
    \item  $\PP \cdot \PP \subset \PP$
    \item  $\LL \cdot \PP \subset \LL$
    \item  $\PP \cdot \UU \subset \UU$
    \item  $\UU \cdot \LL \subset \PP$
    \item  $G$ is the disjoint union $G = \PP \sqcup \PP^{-1} \sqcup \UU \sqcup \LL \sqcup \{e\}$.
\end{enumerate}
\end{prop}

We remark that the disjointness of the union of item (6) would follow from items (1) - (5) if we added the assumption that neither $\UU$ nor $\LL$ contains the identity.  The proof of this proposition is similar to the proof for ordinary left-invariant orders.  If $G$ has a simply connected extension of a partial order, then we make the definitions:
$$\PP = \{g \mid e<g\}, \qquad \UU = \{g \mid e \sim_u g\}, \qquad \LL = \{g \mid e \sim_l g\}.$$
Conversely, if such sets $\PP, \LL$ and $\UU$ exist then we define:
$$g<h  \;\; if \;\;  g^{-1}h \in \PP, \quad  g\sim_u h \;\; if \;\;  g^{-1}h \in \UU, \quad g \sim_l h \;\; if \;\; g^{-1}h \in \LL.$$
Condition (5) is essentially acyclicity and (6) is strong connectivity.

Just as in the cases of other kinds of partial orders, many properties of simply connected extensions may be stated in terms of these sets. For example, the partial order on $G$ given by $P$ is simply connected and needs no extension whenever $G = \PP \cup \PP^{-1} \cup (\PP \cdot \PP^{-1}) \cup (\PP^{-1} \cdot \PP)$.

\section{Groups with Partial Orders act on Order Trees}

In this section we consider groups with partial orders which have left-invariant rectifiable simply connected extensions. We will prove:

\begin{thm}\label{T:order.to.action}
If $G$ is a countable group which admits a partial order with a nontrivial left-invariant rectifiable simply connected extension, then $G$ acts faithfully, nontrivially, and without fixing a unique end by orientation preserving homeomorphisms on an oriented order tree, and hence on a simply connected 1-manifold, $M$. Moreover, $M$ can be chosen to have a point $x$ with trivial stabilizer, so that identifying $G$ with the orbit $Gx$ induces the original partial
order on $G$.
\end{thm}

We will prove this theorem by constructing an oriented order tree $T$ on which the group acts by orientation preserving homeomorphisms. In the construction, a subset of $T$ will be labelled by group elements, and  distinct elements of the group will label distinct points of $T$, so the action will be faithful. We model our construction on the proof of Theorem 6.8 of \cite{Gh}, in which points on the real line are labelled by the elements of $G$ in order to construct the desired action. In our construction, we will label points on intervals of the real line by group elements various in between sets. In order to construct on oriented tree, we think of the equivalence classes in between sets as having orientations induced by the partial order, but maintaining this sense is confusing in the case of an equivalence class consisting of a single group element. We can solve this problem by simply blowing up the set $G$ to a larger poset, $G^+$.  We now define this poset and record some of its basic properties.

\begin{definition}
Let $G^+= \{g,g_-,g_+|g \in G\}$
\end{definition}
\begin{definition}($G$ action on $G^+$)
If $g,h \in G$, define an action of $G$ on $G^+$ via $g(h)=gh, g(h_-)=(gh)_-, g(h_+)=(gh)_+$.
\end{definition}
\begin{definition}\label{D:G+order}(Order on $G^+$)$\forall x \in G$,  $x_-<x<x_+$, and if $x<y$ then $x_+<y_-$ (this is like a Denjoy blow-up done at each of the countably many points).
\end{definition}

It is not hard to see that Definition~\ref{D:G+order} gives $G^+$ a left-invariant partial order.  This new order can be extended following the manner in which the order on $G$ was extended.  To define this extension, let $g_{\epsilon},h_{\delta} \in G^+$ where $\epsilon$ and $\delta$ are $+,-$, or absent. If $g_{\epsilon}\sim h_{\delta}$, then $g \sim h$ in $G$. Furthermore, $g_{\epsilon}$ and $h_{\delta}$ share an upper or lower bound in $G^+$ if and only if $g$ and $h$ share one in $G$.

\begin{definition}
For $h_\delta$ and $g_\epsilon$ as above, define a formal extension of the order on $G^+$ by declaring $g_{\epsilon}\sim_u h_{\delta}$ if $g\sim_u h$ in the extension of the order on $G$ and $g_{\epsilon}\sim_l h_{\delta}$ if $g\sim_l h$ in the extension of the order on $G$.
\end{definition}

\begin{lemma} If $G$ is a group with a rectifiable simply connected extension of a partial order, then the natural order on $G^+$ also is a rectifiable simply connected extension of the corresponding order on $G^+$. In addition, the left $G$-action preserves the extended partial order and no $O_{x,y}$ class consists of a single point as long as $x \neq y$ are both in $G$.
\end{lemma}

\begin{proof}
 The only change in structure is that any time $b$ appeared in, say, a set $B_{x,y} \subset G$, the triple $b_-<b<b_+$ appears in the
 corresponding set $B_{x,y} \subset G^+$. So in particular, the number of equivalence classes under $O_{x,y}$ does not change, and there are no equivalence
 classes consisting of only one point, except possibly at the endpoints if $x$ or $y \in G^+ - G$. As previously mentioned, the $G$-action clearly preserves the order.
\end{proof}

\begin{lemma}\label{l:augmentedbetween} Let $a,b \in G$. Then
\begin{eqnarray*}
a<b &\iff& \{a_+,b_-\} \subset B_{a_-,b_+} \\
a \sim_u b &\iff &\{a_+,b_+\} \subset B_{a_-,b_-}\\
a \sim_l b &\iff &\{a_-,b_-\} \subset B_{a_+,b_+}
\end{eqnarray*}
\end{lemma}

\begin{definition} Let $R$ be the relation on the augmented group $G^+$ defined by: $xRy$ if $B_{x,y}=\{x, y\}$, where $x,y \in G^+-G$, and $xRx$ if $x \in G$.
\end{definition}

\begin{lemma}
R is an equivalence relation.
\end{lemma}

\begin{proof}
$R$ is clearly reflexive and symmetric.  So the only nontrivial case to check for transitivity is if we have $g_sRh_t$ and $h_tRf_r$ where $g,h,f \in G$ and $s,t,r \in \{+,-\}$ and $g_s \neq h_t, h_t \neq f_r$. Now $B_{g_s,f_r}\subseteq B_{g_s,h_t} \cup B_{h_t,f_r}= \{g_s,h_t,f_r \}$, we claim the only possibility is for $g_s=f_r$, hence $g_sRf_r$, which implies that $R$ is transitive. To see this, if $g \neq f$, then by Lemma~\ref{l:augmentedbetween}, $B_{g_s,f_r}$ contains at least four elements of $G^+$, impossible as we saw above that it can contain at most three elements. Hence it must be that $f=g$, but if also $r \neq s$, then $G_{g_s,f_r}=B_{g_s,g_r}=\{g_s, g, g_r \}$, also impossible since we saw above that it contained no elements of $G$.
\end{proof}

We now construct an oriented order tree on which the group acts. We follow the ideas in section 3 of \cite{GK1}, in which order trees are described as countable increasing unions of segment which intersect in a very restricted way.  We begin by constructing a decomposition of the group into a countable union of subsets which will guide the construction of the tree.
If the group $G$ is countable, then $G$ is a countable union of sets $B_{x_i,y_i}=B_i$, where $x_i, y_i \in G$. Consider $x_1$. For each index $i$, include $B_{x_1,x_i}$ in the countable union, and re-index to put each of these right before the corresponding $B_{x_i,y_i}$. Then after re-indexing, we have:

$$G=\bigcup_{i=1}^\infty B_i,$$
$$G_n \cap B_{n+1}\neq \emptyset, \text{where} ~~~G_n=\bigcup_{i=1}^n B_i,$$
$$\forall a,b \in G_n, B_{a,b} \subseteq G_n.$$

\vspace{.2 in} The instructions for the construction will be encoded in the intersections $G_n \cap B_{n+1}$, so we characterize these intersections. Notice that if both $x_{n+1}$ and $y_{n+1}$ are in $G_n$, then $G_{n}=G_{n+1}$ since the entire set $B_{n+1} \subset G_n$.  So after eliminating such indices we can re-index so that $G_n \cap B_{n+1}$ contains $x_{n+1}$ and is an initial segment of $B_{n+1}$ in the ordering $\preceq$ of the set $B_{n+1}$ with $x_{n+1} \prec y_{n+1}$. So the only possibilities for the set $G_n \cap B_{n+1}$ are:

\begin{enumerate}
\item $G_n \cap B_{n+1}= \{x_{n+1}\}$
\item  $G_n \cap B_{n+1}$ is half open, i.e. has more than one element, but is not of the form (3)
\item $G_n \cap B_{n+1}= B_{x_{n+1},z_{n+1}}$ for some $z_n \in G - G_n$.
\end{enumerate}

Now if possibility (3) is the case, just replace $B_{n+1}$ by $B_{z_{n+1},y_{n+1}}$, which reduces the form of the intersection to case (1). So we may assume that each intersection is of type (1) or type (2). Furthermore, if (2) is the case, the furthest equivalence class from $x_{n+1}$ along $B_{n+1}$ which contains elements from $G_n$ has no greatest element in $G_n$, so if we choose one such element, say $z_{n+1}$, and then replace $B_{n+1}$ by $B_{z_{n+1},y_{n+1}}$ as in case (3) above and rename, we now have reduced to the following two cases for the form of $G_n \cap B_{n+1}$:

\begin{enumerate}
\item $G_n \cap B_{n+1}=\{x_{n+1}\}$
\item $G_n \cap B_{n+1}$ is half open and totally ordered in $G_n$.
\end{enumerate}

We will now use this expression of $G$ as the union of the sets $B_i$ to construct a tree $T$ with a labelling by $\nu: G^+ \to T$ as follows. $T$ is the quotient space of a disjoint union $\bigcup_{i=0}^\infty I_i$ modulo a set of identifications $\{R_i\}_{i=1}^\infty$. Each $I_i$ is a closed compact subinterval of the real line, and the identification $R_i$ identifies one endpoint of $I_i$ with one point in the disjoint union $\bigcup_{k=0}^{i-1}I_k$.  We define $T_n$ as the quotient space $\bigcup_{i=0}^n I_i / \{R_1, \ldots, R_n\}$, so that $T_{n+1}=(T_n \cup I_{n+1})/R_{n+1}$.  Then for any $n$, $T= (T_n \cup \bigcup_{i=n+1}^\infty I_i)/\{R_i\}_{i=n+1}^\infty$.  Thus, to define $T$, it suffices to construct $T_0$ and then inductively construct $T_{n+1}$ from $T_n$.
 
To define the labelling function $\nu$, we let $G_n^+=\{g_+,g,g_- |g \in G_n\}$, and we will label a subset of points in $I_n$ by $G_n^+-G_{n-1}^+$.  Apart from a technical detail, this will be the labelling $\tilde\nu_n$ defined below.  Since $G^+$ is the disjoint union of all these subsets, we let $\nu$ be the labelling of $T$ by $G^+$ obtained by defining $\nu$ restricted to each disjoint subset $G_n^+-G_{n-1}^+$ to be the appropriate labelling map above, and we define $\nu_n$ be the labelling of $T_n$ obtained similarly. Finally, each tree $T_n$ is naturally homeomorphic to a subtree of both $T_{n+i}$ and $T$, and we will sometimes abuse notation and refer to all of these homeomorphic copies by the same name, $T_n$.

\begin{definition}(The Construction)
Let $G$ be a group with a left-invariant rectifiable simply
connected extension of a partial order. Express the group as a
union of the form:
$$G=\bigcup_{i=1}^\infty B_i,~~\text{where}~~B_i=B_{x_i,y_i}.$$
Further assume that,
$$G_n \cap B_{n+1}\neq \emptyset, ~~~\text{where}~~G_n=\bigcup_{i=1}^n B_i,$$
$$\forall a,b \in G_n, B_{a,b} \subseteq G_n,$$
and in addition $G_n \cap B_{n+1}$ contains $x_{n+1}$ and has one of the two following forms:

\begin{enumerate}
\item $G_n \cap B_{n+1}=\{x_{n+1}\}$
\item $G_n \cap B_{n+1}$ is half open and totally ordered in $G_n$.
\end{enumerate}

Base step: Let $T_0=I_0$ be a copy of the unit interval, and define $\nu_0: \{(x_1)_-,x_1,(x_1)_+\} \rightarrow T_0$ by $\nu_0((x_1)_-)=0$, $\nu_0(x_1)=0.5$, $ \nu_0((x_1)_+)=1$.

Inductive step: Suppose $T_n$ has been constructed, along with $\nu_n:G_n^+ \rightarrow T_n$. We must label a closed segment of the real line, $I_{n+1}$ by the elements of $G_{n+1}^+-G_n^+$ and then specify the identification $R_{n+1}$ to construct $T_{n+1}$. The new labelling $\nu_{n+1}:G_{n+1}^+ \rightarrow T_{n+1}$ will be defined to be the map $\nu_n$ restricted to $G_n^+$ and the labelling specified below of $I_{n+1}$ restricted to $G_{n+1}^+-G_n^+$. The details will depend upon the form of $G_n \cap B_{n+1}$.

Case 1: $B_{n+1} \cap G_n= \{x_{n+1} \}$.  To simplify notation, suppress subscripts and abbreviate $x_{n+1}$ by $x$ and $y_{n+1}$ by $y$. Now suppose $B_{n+1}$ is made up of $k$ equivalence classes under $O_{x,y}$; each is totally ordered with respect to the partial order on $G$, and the classes themselves have a natural linear order along $B_{x,y}$. We index these equivalence classes by $i$ where $1 \leq i \leq k$, where we take the first class, or $i=1$, to be the one containing $x$. We construct a labelling of a copy of $[0,k] \subset \mathbb{R}$. For each index $i$, augment the equivalence class by adding, for each $g$ in the class, $g_-$ and $g_+$. First, identify pairs of elements which are equivalent mod $R$. The resulting set is still totally ordered.  Now define $(\nu_{n+1})_i$ from this totally ordered set of labels to $[i-1,i]$ as in the classical construction in the proof of Theorem 6.8 of \cite{Gh}, which lays a totally ordered set down on the real line. However, if this modified equivalence class happens to have least (or greatest) elements, label $i-1$ (or $i$ respectively) by these elements. If not, the endpoints $i$ and $i+1$ receive no labels, but will be limit points of labelled points. Notice that since no augmented equivalence class has only one element, there is never confusion about whether to label $i-1$ or $i$ by a given extremal element of the modified equivalence class. Since $B_{n+1} \cap G_n= \{x \}$, if we take $\tilde{\nu}_{n+1}$ to be the union of all these maps, the domain of $\tilde{\nu}_{n+1}$ intersects the domain of $\nu_n$ only at the triple $x_-, x, x_+$.  Exactly one of these augmented elements $x_-,x_+$ is between $x$ and $y$, for ease of discussion suppose it is $x_-$. Then, on the interval $[0,1]$, $0=\tilde{\nu}_{n+1}(x_+)<\tilde{\nu}_{n+1}(x)<\tilde{\nu}_{n+1}(x_-) \leq 1$,  and furthermore, the open segments $(0,\tilde{\nu}_{n+1}(x))$ and $(\tilde{\nu}_{n+1}(x),\tilde{\nu}_{n+1}(x_-))$ contain no images under $\tilde{\nu}_{n+1}$. We remove the segment $[0,\tilde{\nu}_{n+1}(x_-))$ from $[0,k]$, leaving just $[\tilde{\nu}_{n+1}(x_-),k]$. This segment, $[\tilde{\nu}_{n+1}(x_-),k]$, is the segment $I_{n+1}$.  It is possible that $\tilde{\nu}_{n+1}(x_-)=1$ in the case that there were no other group elements besides $x$ in the equivalence class of $x$ along $B_{n+1}$. Otherwise, $\tilde{\nu}_{n+1}(x_-)<1$. 

To define the gluing relation $R_{n+1}$, identify the point $\tilde{\nu}_{n+1}(x_-)$ to the point $\nu_n(x_-)$ in $T_n$ to form $T_{n+1}$, i.e. define $T_{n+1}= T_n \cup [\tilde{\nu}_{n+1}(x_-),k] /\tilde{ \nu}_{n+1}(x_-)=\nu_n(x_-)$, and define $\nu_{n+1}:G_{n+1}^+ \rightarrow T_{n+1}$ as just $\nu_n$  restricted to $G_n^+$ and $\tilde{\nu}_{n+1}$ restricted to $G_{n+1}^+-G_n^+$, as in Figure~\ref{F:tnplus1.1}. Notice that in fact, the domains of $\tilde{\nu}_{n+1}$ and $\nu_n$ intersect in precisely the element $x_-$, and the relation $R_{n+1}$ identifies $\tilde{\nu}_{n+1}(x_-)$ with $\nu_n(x_-)$


\begin{figure}
\begin{center}
\psfrag{a}{\small $\nu_n(x_+)$}
\psfrag{b}{\small $\nu_n(x)$}
\psfrag{c}{\small $\nu_{n+1}(x_-)=0$}
\psfrag{1}{\small $1$}
\psfrag{2}{\small $2$}
\psfrag{l}{\small $k-1$}
\psfrag{k}{\small $k = \nu_{n+1}(y_\epsilon)$}
\psfrag{G}{\small $T_{n}$}
\psfrag{B}{\small $B_{n+1}$}
\psfrag{cdots}{\small $\cdots$}
    \includegraphics[width = 3in]{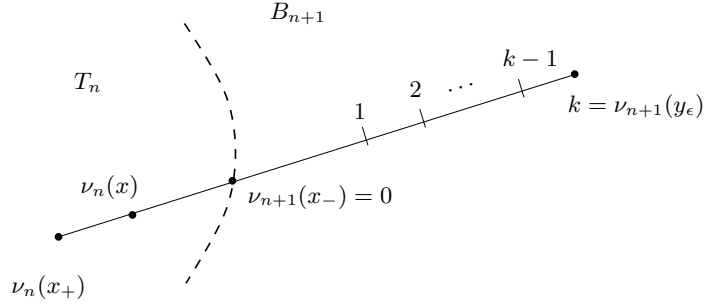}
    \caption{$T_{n+1}$ after the gluing of $B_{n+1}$, Case 1}\label{F:tnplus1.1}
\end{center}
\end{figure}

Case 2: $B_{n+1} \cap G_n$ is half open and totally ordered in $G$. So $B_{n+1} \cap G_n$ is a totally ordered set of group elements, and $x$ is either the greatest or the least of them. For ease of discussion, assume it is the least, so that there is no greatest element. Then choose a countable, increasing subsequence $\{a_i\} \subset B_{n+1} \cap G_n$ which is not bounded above (with respect to the total order) by any element in $B_{n+1} \cap G_n$. So these group elements label points in $T_n$, and since $T_n$ is a finite union of compact intervals, there is some limit point $l$ of the sequence $\{\nu_n(a_i)\}$ in $T_n$.  Now consider the set $B_{n+1} - (G_n \cap B_{n+1})$, which intersects only finitely many (say $k$) of the equivalence classes in $B_{n+1}$.  Label a copy of $[0,k]$ by elements in $G_{n+1}^+-G_n^+$ as in Case 1, with either $\tilde{\nu}_{n+1}(y_-)=k$ or $\tilde{\nu}_{n+1}(y_+)=k$.  In this case, the segment $[0,k]$ is $I_{n+1}$.  If an element $g \in B_{n+1} - (G_n \cap B_{n+1})$ closest to $x$ exists, then let $\tilde{\nu}_{n+1}(g_s)=0$, $s \in \{+,-\}$, where $g_s \in B_{x,g}$ in the augmented group. If no such $g$ exists, $0$ remains unlabelled. For every other $h$ in that equivalence class, points in $(0,1]$ are labelled by $h_-,h,h_+$ as usual. Finally, define the gluing relation $R_{n+1}$ by identifying the point $0\in I_{n+1}$ to the point $l\in T_n$, i.e. define $T_{n+1}= T_n \cup I_{n+1} / l=0$, and define $\nu_{n+1}:G_{n+1}^+ \rightarrow T_{n+1}$ to agree with  $\nu_n$ when restricted to $G_n^+$ and $\tilde{\nu}_{n+1}$ restricted to $G_{n+1}^+-G_n^+$, as in Figure~\ref{F:tnplus1.2}.


\begin{figure}
\begin{center}
\psfrag{a}{\small $\nu_n(x_+)$}
\psfrag{b}{\small $\nu_n(x)$}
\psfrag{c}{\small $l=0$}
\psfrag{a1}{\small $\nu_n(a_1)$}
\psfrag{a2}{\small $\nu_n(a_2)$}
\psfrag{1}{\small $1$}
\psfrag{2}{\small $2$}
\psfrag{l}{\small $k-1$}
\psfrag{k}{\small $k = \nu_{n+1}(y_\epsilon)$}
\psfrag{G}{\small $T_{n}$}
\psfrag{B}{\small $B_{n+1}$}
\psfrag{cdots}{\small $\cdots$}
    \includegraphics[width = 3in]{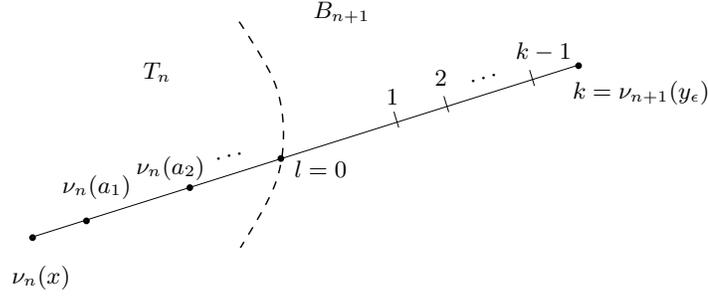}
    \caption{$T_{n+1}$ after the gluing of $B_{n+1}$, Case 2}\label{F:tnplus1.2}
\end{center}
\end{figure}

\end{definition}

\begin{figure}
\begin{center}
\psfrag{a}{\small $\nu_n(a)$}
\psfrag{b}{\small $\nu_n(d_-)$}
\psfrag{c}{\small $\nu_{n}(d)$}
\psfrag{d}{\small $\nu_n(c_+) = \nu_n(d_+)$}
\psfrag{e}{\small $\nu_n(c)$}
\psfrag{f}{\small $\nu_n(c_-)$}
\psfrag{g}{\small $\nu_{n}(b)$}
    \includegraphics[width = 3in]{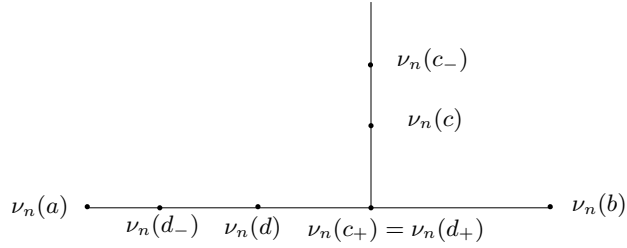}
    \caption{One possibility for multiple labellings in Property 4}\label{F:property4}
\end{center}
\end{figure}

One shows inductively that $T_n$ has the following properties:

\begin{lemma} In the construction above, for each $n$ we have:

\begin{enumerate}

\item $T_n$ is a tree.

\item The connected components of the complement of the closure of $\nu_n(G_n^+)$ are precisely the intervals of the form $(\nu_n(g),\nu_n(g_-))$ or $(\nu_n(g),\nu_n(g_+))$, for some $g \in G_n$.

\item $\forall a,b \in G_n$, the unique interval $[\nu_n(a),\nu_n(b)]$ in $T_n$ contains points, labelled in the natural order by the elements of the set $B_{a,b} \subset G_n^+$. Furthermore, no other elements of $G$ can label these points, and if another element $c \in (G_n^+-G_n) - B_{a,b}$ labels a point along $[\nu_n(a),\nu_n(b)]$ both of the following conditions are satisfied:
    \begin{enumerate}
    \item  Either $cRd$ for $d \in B_{a,c}\cap B_{a,b} \subset G_n^+$, or $\nu_n(c)$ is a limit point of points labelled by elements of $B_{a,c}\cap B_{a,b}$.
    \item  Either $cRd$ for $d \in B_{b,c}\cap B_{a,b} \subset G_n^+$, or $\nu_n(c)$ is a limit point of points labelled by elements of $B_{b,c}\cap B_{a,b}$.  
    \end{enumerate}
An example of the first possibility occurring on the $\nu_n(a)$ side is illustrated in Figure~\ref{F:property4}
\item For $x,y \in G_n^+$, $\nu_n(x)=\nu_n(y)\iff xRy$.

\end{enumerate}
\end{lemma}

Note that the first property shows that the space $T$ is simply connected.  The fourth will ensure that the labelling by $G^+$, though not injective, induces an injective map from $G^+/R$ to $T$, so that in particular distinct group elements label distinct points of $T$. The second and third properties will ensure that the natural action of $G$ on the set of labelled points extends to an orientation preserving action on the space $T$.


We now prove,

\begin{prop}\label{P:thetree} The space $T$ constructed above has the structure of an oriented order tree, and the group $G$ acts faithfully, nontrivially, without fixing a unique end, by orientation preserving homeomorphisms on $T$.
\end{prop}
\begin{proof} To see that $T$ has the structure of an oriented order tree, consider the set of the images of the closed segments $[j,j+1]\subset I_i$ for all integers $i,j \geq 0$.  Recall that in Case 1 of the construction, the interval $[0,k]$ was shortened by removing an open subinterval and leaving, $I_i = [\tilde{\nu}_{i}(x_-),k]$.  Therefore, the initial subinterval of $I_i$ that we wish to consider in Case 1 is $[\tilde{\nu}_{i}(x_-),1]$, rather than $[0,1]$.  For the rest of the proof, we will abuse notation by referring to the left endpoint of any such interval, $[\tilde{\nu}_{i}(x_-),1]$, by 0. Denote by $\sigma$ the image in $T$ of one of the above closed intervals $[j,j+1]$. By construction, each of these has at least one interior point labelled by a group element $g$. If $g_-$ labels a point in  $[j,\nu(g))$ let $i(\sigma)$ be the endpoint $j$ and let $f(\sigma)$ be $j+1$. On the other hand, if $g_-$ labels a point in the image of $(\nu(g),j+1]$, let $i(\sigma)$ be $j+1$ and let $f(\sigma)$ be $j$. Note that this definition is independent of the choice of $g$, since the group elements which are labels in $[j,j+1]$ are totally ordered. Let $S^+$ be the set of all such $\sigma$, together with their closed subsegments and unions of any pairs $\{\sigma_1,\sigma_2\}$ where $\sigma_1\cap \sigma_2=f(\sigma_1)=i(\sigma_2)$ (with the obvious assignments of initial and final endpoints). This set of positively ordered segments gives $T$ the structure of an oriented order tree.

To establish the action we will use the labelling $\nu:G^+ \rightarrow T$.  Note that $\nu(g)=\nu(h) \Rightarrow g=h$ for $g,h \in G$, and that the set $\nu(G)$ spans $T$. We term a branch point any point in the tree which was the image of an integer point or an endpoint of $I_n$ for some $n$. Branch points, by construction, are always either labelled by elements of $G^+-G$ or are limit points of labelled points. Note that this set of branch points includes all points where the tree genuinely branches, in addition to any point where the orientation changes. Hence if a connected subset of $T$ contains no branch points, it is homeomorphic by order preserving homeomorphism to an interval of the real line.

To see that the group acts on $T$, let $L$ be the subset of $T$ which is labelled by elements of $G^+$. Consider the complement of the closure of $L$. Since the complement of the closure of $L$ contains no branch points, any connected component of this complement is an open interval of the form $(p,q)$. We claim that $p$ and $q$ are in $L$. For suppose that $x \in (p,q)$. Then $x \in [\nu(g),\nu(h)]$ for some $g,h \in G$. Now $g,h \in G_n$ for some $n$, hence $[\nu (g),\nu (h)] \in T_n$. But since $x$ is unlabelled, it is certainly unlabelled by $\nu_n$ in $T_n$, so $x \in (\nu_n(f),\nu_n(f_s))$ for some $f \in G_n, s \in \{+,-\}$. But note that in $T_{n+i}$, points in the segment $(\nu_{n+i}(f),\nu_{n+i}(f_s))$ can never be labelled, for by property 3 of the construction, they could only be labelled by elements in $G^+-G_n^+$ if they were limit points of labelled points, or if they were already labelled. Hence, since $\nu(f), \nu(f_s) \in L$, $(p,q)=(\nu(f),\nu(f_s))$.

Now since $G$ acts on $G^+$ respecting the relation $R$, $G$ clearly acts on $L$, and we would like to extend this action continuously to the closure of $L$. Let $x$ be a point in this closure but not in $L$. Then we may have many segments of the form $[a,x]$, any two of which intersect only at $x$, each containing an increasing (in the total order on $[a,x]$ with $x$ being greatest) sequence of labelled points converging to $x$. Suppose $[a,x]$ and $[b,x]$ are any two such segments, and let $\{\nu(a_i)\}$ and $\{\nu(b_i)\}$ be the two sequences of points. Note that for some $n$, $[a,x]\cup [x,b] \in T_n$. Now let $g$ be an arbitrary group element, and we worry about ambiguity in defining $gx$. The worry is that after transforming by the group element $g$, the two sets of points $\{\nu(ga_i)\}$ and $\{\nu(gb_i)\}$  along $[\nu(ga),\nu(gb)]$ will have different limit points, say $x_a$ and $x_b$. We claim there can be no labelled points in the open segment $(x_a,x_b)$. If there were, then that label (say $h$) would lie in the set $B_{ga_i,gb_j}\subset G^+$ for every $i$ and $j$. But applying $g^{-1}$ we see that $g^{-1}h \in B_{a_i,b_j}$ for every $i$ and $j$. But then $g^{-1}h$ must label a point along $[a,b]$ between the points $\nu(a_i)$ and $\nu(b_j)$ for every $i$ and $j$. However, $x$ is the only such point, and it is unlabelled by hypothesis. So $x_a$ and $x_b$ must be endpoints of some connected component of the complement of the closure of $L$, hence they are actually labelled points. Since they are different points, the labels are distinct, and not even equivalent modulo $R$. But then applying $g^{-1}$ again, we produce distinct labels not equivalent modulo R which are between $a_i$ and $b_j$ for any $i$ and $j$, hence they label distinct points in $[\nu(a_i),\nu(b_j)]$ for every $i$ and $j$. But this is nonsense, as there is only one point between $[\nu(a_i),\nu(b_j)]$ for every $i$ and $j$. So the group action may be extended continuously to the closure of $L$. Since the complement of the closure is a disjoint union of open intervals, extend the action to all of $T$ by extending linearly across these intervals.
The action clearly preserves the orientation,  as $g(h_-)=(gh)_-$ and $g(h_+)=(gh)_+$ for $g, h \in G$.

Since the labelling $\nu$ restricted to $G$ is injective, the action is faithful. The fact that it is has neither a global fixed point nor a unique fixed end follows from the nontriviality of the extension of the partial order.
First, suppose to the contrary that $x \in T$ is fixed by all $g \in G$. First, $x$ cannot be in the complement of the closure of $L$, for if it were, then $x \in (\nu(h_s), \nu(h))$ for some $h \in G$, which implies that $\forall g \in G, gx \in (\nu((gh)_s),\nu(gh))$. So $gx=x \Rightarrow gh=h \Rightarrow g=e$. But $x$ cannot be in $L$ itself, for if $x \in L$, then $x=\nu(h_s)$ for some $h \in G$, $s \in \{+,-\}$. But then $\forall g \in G, \nu(h_s)=g(\nu(h_s))=\nu((gh)_s)$. So $\forall g \in G, (gh)_sRh_s \Rightarrow \forall f,g \in G, f_sRg_s$. Hence either $\forall~~ f,g \in G,~~ f\sim_ug$, or $\forall~~ f,g \in G, ~~f \sim_l g$, or in other words, the extension was trivial. So if $x$ is fixed by all of $G$, it must be that $x$ is a limit point of labelled points, but is itself unlabelled. Then choose any group element $g \in G$, and consider the segment $[\nu(g),x]$. If $[\nu(g),x]$ is not an oriented segment of $T$, because there are at most finitely many switches in orientation along it , you may re-choose $g$ so that $[\nu(g),x]$ is oriented. But given any other $h \in G$, since $g(h^{-1})$ fixes $x$ and preserves orientation, $[\nu(h),x]$ must also be an oriented segment $[\nu(h),x]$ and $[\nu(g),x]$ are either both oriented towards or away from $x$. Hence all segments of the form $[\nu(h),x]$ must be either oriented towards $x$ or oriented away. But then given $f,h \in G$, since $[\nu(f), \nu(h)] \subset [\nu(f),x] \cup [\nu(h), x]$, and since labels along these segments correspond to between sets in $G$, either $f$ and $h$ are comparable, or $f\sim_u h$ in the case $x>\nu(g)~\forall g\in G$ or $f \sim_l h$ in the case that $x<\nu(g)~\forall g\in G$. Note that since $x$ is a limit point of labelled points, by construction there must be some pair $g,h$ where $x \in [\nu(g), \nu(h)]$, so that $g\sim h$, and hence the extension is trivial.

Now suppose the action has a unique fixed end $e$.  For an arbitrary element $g$,  choose a ray $\rho$ from $\nu(g)$ representing $e$.   If $\nu(g_-) \in \rho$ we say that $\nu(g)$ points towards the end $e$ represented by $\rho$, and if $\nu(g_+) \in \rho$ we say $\nu(g)$ points away from $e$. Since the group acts transitively on itself but fixes the end $e$, for any group element $h$, $\nu(h)$ will point the same way with respect to $e$ as $\nu(g)$.  If both $\nu(g)$ and $\nu(h)$ point towards $e$, either $g$ and $h$ are comparable, or $g \sim_l h$; and if both $\nu(g)$ and $\nu(h)$ point away from $e$, they are either comparable, or $g \sim_u h$.  So all incomparable pairs must be of the same type. Notice that there is at least one such pair, for if not, $G$ is totally ordered, and then $T=\mathbb{R}$, and $e$ is not the unique fixed end. So the extension is in fact trivial if there is a unique fixed end.

Properties 3 and 4 also allow us to translate the characterization in Lemma~\ref{l:augmentedbetween} of the order relationships in $G$ in terms of between sets in $G^+$ to a characterization in terms of labels on shortest paths in $T$. \end{proof}

\begin{rmk}Branch points were never labelled by group elements in this construction, so all group elements are in the interiors of some oriented segment.  As usual, one can blow this up to an oriented Non-Hausdorff 1-manifold. Since group elements don't label branch points, the blowing up retains the property that each group element labels a distinct point, and the group will act as desired on the oriented simply connected 1-manifold.
\end{rmk}

\section{Groups acting on oriented order trees have simply connected partial orders}

The goal of this section is to prove a collection of partial converses to Theorem~\ref{T:order.to.action} by showing that groups that act on oriented order trees admit left-invariant partial orders with nontrivial rectifiable simply connected extensions.  We continue to restrict ourselves to $\mathbb{R}$-order trees, but with some technical improvements Theorems~\ref{T:arbitrarytree} and~\ref{T:weakstab} can be extended to arbitrary order trees. This is achieved by replacing ``simply connected oriented 1-manifold" by ``oriented order tree without branching" in Theorems~\ref{T:groupgetsorder} and ~\ref{T:nontrivialstab} and defining the geodesic spine $\gs[x]y$ to be the intersection of all paths from $x$ to $y$. The technical improvements involve dealing with the facts that in an arbitrary order tree, a finite ray need not have any endpoint and a geodesic spine may not decompose into a disjoint union of segments. Indeed, the definition of finite ray would be changed to include any ray that is contained in some finite geodesic spine.

We first remark that there exist many examples of countable groups that act nontrivially, faithfully and without fixing a unique end by order preserving homeomorphisms on an oriented order tree but which cannot admit a left-invariant partial order with nontrivial rectifiable simply connected extension. Therefore, we cannot hope for a direct converse to Theorem~\ref{T:order.to.action}.

To construct an example violating the converse of Theorem~\ref{T:order.to.action}, suppose that $G_1 = \langle g_1 \rangle$ and $G_2 = \langle g_2 \rangle$ are non-isomorphic cyclic groups each containing a proper subgroup isomorphic to $H$.  Form the free product of $G_1$ and $G_2$ amalgamated along $H$; set $G = G_1 *_{H} G_2$.  Standard Bass-Serre theory yields a simplicial tree $T$ with a $G$ action.  The quotient of $T$ by $G$ is a graph with two vertices, $v_1$ and $v_2$, and one edge, $e$.  If we orient $e$ arbitrarily and lift the orientation to a $G$-equivariant orientation of $T$, then $T$ has the structure of an oriented order tree and the action of $G$ is by orientation-preserving homeomorphisms.  Moreover, the action is faithful, nontrivial and fixes no unique end.  However, $G$ cannot admit a partial order with nontrivial rectifiable simply connected extension, which can be seen as follows.

Suppose that $G$ does admit a partial order with nontrivial rectifiable simply connected extension.  Theorem~\ref{T:order.to.action} then applies to $G$, so let $T'$ be the tree constructed in that theorem.  Recall from~\cite{RSS} that if an element $g\in G$ has no fixed point in $T'$ then $g$ acts as translation along some axis, so no nonzero power of $g$ has a fixed point.  Since $g_1$ and $g_2$ are torsion elements in $G$, we may choose fixed points $x$ and $y$ of $g_1$ and $g_2$ respectively.  The points $x$ and $y$ are fixed by all of $G_1$ and $G_2$, respectively, because $g_i$ generates $G_i$. Since $H \subset G_1 \cap G_2$, $H$ must stabilize the geodesic spine $GS_{(x,y)}$.  We will show that the action on the interior of this geodesic spine is faithful.  Since the interior of $GS_{(x,y)}$ is   homeomorphic to $\mathbb R$, this will show that $H$ has a faithful action on $\mathbb R$.  Such an action is impossible since it would induce on $H$ a total order, as in Theorem 6.8 of \cite{Gh}, but a finite group cannot admit a total order.  Therefore, $G$ cannot admit a partial order with nontrivial rectifiable simply connected extension.

To show that the action of $H$ on the interior of $GS_{(x,y)}$ is faithful, it will suffice to show that the interior of $GS_{(x,y)}$ contains a point labelled by an element of $G$, since the stabilizer of such a point is trivial.  Note first that since both endpoints have nontrivial stabilizers, neither $x$ nor $y$ can be labelled by an element of $G$, so if one is labelled, the label must belong to $G^+ - G$.  On the other hand let,
$$A := \bigcup_{g\in G} gGS_{(x,y)}.$$
Now, any $g$ has a decomposition $g=h_1k_1h_2k_2\cdots h_mk_m$ with $h_i \in G_1$ and $k_i \in G_2$, which produces a path,
$$h_1GS_{(x,y)}, \,h_1k_1GS_{(x,y)}, \,h_1k_1h_2GS_{(x,y)} , \ldots ,\,h_1k_1\cdots h_mk_mGS_{(x,y)}$$
from $GS_{(x,y)}$ to $gGS_{(x,y)}$.  Therefore $A$ is a connected $G$-invariant subset of $T'$, so it is an invariant subtree.  Since the action is minimal, $A$ must be all of $T'$.  But, $T'$ contains points labelled by elements in $G$, and this set of points is invariant under the action of $G$.  Since neither $x$ nor $y$ is labelled by an element of $G$, none of their translates are.  Hence, $G$-labelled points must lie in the translates of the interior of $GS_{(x,y)}$, and so a $G$-labelled point must lie in the interior of $GS_{(x,y)}$ itself.

In this example, the groups $G_1$ and $G_2$ stabilize the points $\widetilde{v}_1$ and $\widetilde{v}_2$, which are lifts of the vertices $v_1$ and $v_2$ to $T$. In fact, every point of $T$ has a nontrivial stabilizer. The above discussion illustrates that the fact that these stabilizers are not left orderable constitutes a major obstruction to defining an extended simply connected partial order on $G$ by using its action on $T$.  As we will prove in the most general of these results, Theorem~\ref{T:weakstab}, no other obstruction exists.  Thus, to prove a partial converse to Theorem~\ref{T:order.to.action}, we must impose the condition that there exist a point in $T$ with left-orderable stabilizer.

In structuring these theorems, we distinguish between groups that act minimally on  oriented simply connected 1-manifolds, which we will show have simply connected partial orders, and groups that act minimally on general oriented order trees, which in general have extensions of simply connected partial orders.

The group theoretic motivation for this distinction can be seen in the example in Section~\ref{s:order} just before Definition~\ref{D:extension}. This example shows an action of the infinite dihedral group on an order tree that gives rise to an acyclic (but not strongly connected) partial order with rectifiable simply connected extension. Moreover, this action induces an action on a simply connected 1-manifold which we have shown is not minimal, and has no minimal invariant submanifold, corresponding to the fact that the infinite dihedral group cannot admit a full simply connected partial order.

 In the case of 3-manifold groups, this distinction has a natural topological interpretation.  Namely, a minimal group action on an oriented simply connected 1-manifold arises from the presence of a minimal foliation in the 3-manifold, whereas a more general minimal action on an oriented order tree induces a (not necessarily minimal) action on an oriented 1-manifold, which arises from the presence of a (not necessarily minimal) foliation in the 3-manifold.

For the case of a minimal action on a 1-manifold (or later on, the case of a minimal action on an order tree), we first consider first an action with a trivial stabilizer in Theorem \ref{T:groupgetsorder} (in Theorem \ref{T:arbitrarytree}, respectively), and then extend to actions with a left orderable stabilizer in Theorem \ref{T:nontrivialstab} (in Theorem~\ref{T:weakstab}, respectively).

For minimal actions on simply connected 1-manifolds, we will prove:

\begin{thm}\label{T:groupgetsorder}
If a countable group $G$ acts minimally, and without fixing a unique end, on an oriented, simply connected 1-manifold $T$ by orientation preserving homeomorphisms, and there is some point in the manifold with trivial stabilizer, then $G$ admits a nontrivial left-invariant rectifiable simply connected partial order.
\end{thm}

\begin{definition}\label{D:pointing}  If $T$ is an oriented non-Hausdorff simply connected $1$-manifold, we say that $x$ {\em points at} $y$ if $y\in x^-$, and $x$ {\em points away from} $y$ if $y\in x^+$.
\end{definition}

Throughout the rest of this section whenever we have a group $G$ acting on an order tree, $T$, we will denote the orbit under $G$ of the element of $x\in T$ by $\OO(x)$.

\begin{lemma}\label{L:lowerbd}
Let $T$ be an oriented non-Hausdorff simply connected $1$-manifold and suppose that $G$ acts minimally on $T$. If $x,y \in \mathcal{O}(v)$ are incomparable and share a lower bound in $T$ then $x$ and $y$ share a lower bound in $\mathcal{O}(v)$, and if $x,y \in \mathcal{O}(v)$ share an upper bound in $T$ then $x$ and $y$ share an upper bound in $\mathcal{O}(v)$
\end{lemma}

\begin{proof}

We consider only the case of lower bounds;  the case for upper bounds is similar.  First, we show that there can be no element $a\in T$ towards which every element of $\mathcal{O}(v)$ points.  Suppose that every element of $\mathcal O(v)$ points at the element $a \in T$.  Then every element of $\mathcal{O}(v)$ points at every element of $\mathcal{O}(a)$, and in particular $a \notin \mathcal O(v)$. Let $A=\bigcup_{\alpha,\beta \in \mathcal{O}(a)} GS_{(\alpha,\beta)}$.  Then $A$ is a G-invariant implicit subtree, but no element of $\mathcal{O}(v)$ can be in $A$, since it would have to point at both ends of some geodesic spine $GS_{(\alpha,\beta)}$.  This contradicts the fact that the action is minimal, so there can be no such element $a$.

Suppose now that $x$ and $y$ are incomparable and share a lower bound in $T$ but share no lower bound in $\mathcal{O}(v)$.  We claim that in this case there is a point $a\in T -\mathcal{O}(v)$ such that every point of $\mathcal{O}(v)$ points at $a$.  By the previous paragraph, this is impossible, which will prove that if $x$ and $y$ have a lower bound in $T$, they have one in $\mathcal{O}(v)$ as well.

Recall that the geodesic spine $GS_{(x,y)}$ contains finitely many cusp pairs.  If an interior segment of $GS_{(x,y)}$ contains a point of $ z\in \mathcal{O}(v)$, $z$ is incomparable to either $x$ or $y$, and shares a lower bound in $T$ with this element.  Suppose that it is $y$ that shares the lower bound with $z$.  Any lower bound of $y$ and $z$ is also a lower bound of $x$ and $z$.   Therefore, by possibly replacing $x$ and $y$ with other elements of $\mathcal{O}(v)$ on the segment of $GS_{(x,y)}$, we may assume that no interior segment of $GS_{(x,y)}$ contains points of $\mathcal{O}(v)$.

Let $t$ be the first point in a cusp pair along $GS_{(x,y)}$ on the way from $x$ to $y$.  Choose an open neighborhood $N'$ of $t$ homeomorphic to $\mathbb{R}$.  Let $N$ be the open segment of $N'$ not contained in $GS_{(x,y)}$.  Then $N$ is directed away from both $x$ and $y$, so $N$ can contain no points of $\mathcal{O}(v)$.  Finally, since no interior segment of $GS_{(x,y)}$ intersects $\mathcal{O}(v)$, every point of $GS_{(x,y)} \cap \mathcal{O}(v)$ points at every point of $N$.

Let $a \in N$.  We will show that everything in $\mathcal{O}(v)$ points at $a$.  Since we already have everything in $\mathcal{O}(v) \cap GS_{(x,y)}$ pointing at $N$, let $z \in \mathcal{O}(v) - GS_{(x,y)}$.  Since $z$ does not lie in $N$ or $\gs[x]y$, the three elements, $x,y$ and $a$ all lie in the some component of $T-\{z\}$.   There are two cases to consider.

Case 1:  $z\in x^-\cap y^-$.  Since $x$ and $y$ share no lower bound in $\mathcal{O}(v)$, $x,y \in z^-$ (otherwise we would have $x>z, y>z$).  Therefore $a$ must also lie in $z^-$, which means that $z$ points at $a$.

Case 2:  $z\not\in x^-\cap y^-$.  We consider the case that $z \in x^+$;  the case $z\in y^+$ is handled similarly.  We wish to show that $z$ points at $x$ and hence at $a$.  Choose $g\in G$ with $w:=gx$ in the component of $T-a$ that is contained in $x^-\cap y^-$.  By Case 1, $a\in w^-$.  And, since $z\in x^+$, we have $gz\in w^+$.  Therefore, $a$ and $gz$ lie in opposite components of $T-w$, so $a$ and $w$ lie in the same component of $T-(gz)$.

Now, $w\in x^-\cap y^-$ so $w^+ \subset x^-\cap y^-$ (otherwise we would have $w^- \subset x^-\cap y^-$ forcing $w$ to be a lower bound for $x$ and $y$.)  Therefore, $gz\in w^+ \subset x^- \cap y^-$ and by Case 1, $a\in (gz)^-$.  Since $a$ and $w$ lie in the same component of $T-(gz)$, we have $w\in (gz)^-$.  In other words, $gx\in (gz)^-$ so that $x\in z^-$.  Since $z\in x^+$ (by the assumption for Case 2), we have $x^-\subset z^-$.  Since $a\in x^-$, we have $a\in z^-$ showing that $z$ points at $a$.
\end{proof}

We now prove Theorem~\ref{T:groupgetsorder}.

\begin{proof}[Proof of Theorem~\ref{T:groupgetsorder}]
Since $T$ is an oriented simply connected $1$-manifold, $T$ is a simply connected poset. Choose $x\in T$ such that $Stab(x)$ is trivial.  By Lemma~\ref{L:lowerbd} the subposet $\mathcal O(x)$ is strongly connected.  Since $\mathcal O(x)$ is a subposet of the acyclic poset $T$, $\mathcal O(x)$ is itself acyclic.  Therefore, the partial order of $\mathcal O(x)$ is simply connected.  Also, since the poset $T$ is rectifiable, so is $\mathcal O(x)$. Therefore, the left-invariant partial order given to $G$ by identifying it with $\mathcal O(x)$ is rectifiable and simply connected.  The rest of the proof is devoted to proving that the simply connected partial order of $\OO(x)$ (and therefore of $G$) is nontrivial.

We will show that the assumption that the order is trivial leads to a contradiction.  This is done by showing that under this assumption, $T$ must be homeomorphic to $\mathbb R$.  Since $G$ acts faithfully on $T$, this shows that the order on $\OO(x)$ would be total, and hence nontrivial.

Assume towards a contradiction that the order on $\OO(x)$ is trivial.  Since a total order is not trivial, $\OO(x)$ must contain non-comparable pairs of elements.  Moreover, all non-comparable pairs must be of the same type; either all such pairs satisfy $\sim_l$ or all satisfy $\sim_u$.    We consider only the case in which each two non-comparable elements $x,y \in \OO(x)$ satisfy $x\sim_l y$ (the case for $x\sim_u y$ is similar).

Our first step towards proving that $T$ is a line is to show that there exists a pair of comparable elements in $\OO(x)$.  Suppose not.  Then, $y\sim_l z$ for any $y,z \in \mathcal{O}(x)$.  Therefore, for any $y,z\in \mathcal{O}(x)$, $\mathcal{O}(x) \cap B_{y,z} = \{y,z\}$.  But $GS_{(y,z)} = B_{y,z}$,  so no element of $\mathcal{O}(x)$ ever separates two other elements of $\mathcal{O}(x)$.  Therefore the set,
$$I := \left( \bigcup_{y,z \in \mathcal{O}(x)} GS_{(y,z)} \right) - \mathcal{O}(x)$$
is a $G$-invariant implicit subtree of $T$.  Since the $G$ acts minimally on $T$, $I = \emptyset$.  Therefore, no two points of $\mathcal{O}(x)$ are separable, so that $\mathcal{O}(x)$ itself is a proper $G$-invariant implicit subtree of $T$, again contradicting the minimality of the action.  Thus, there must be at least two comparable elements of $\mathcal{O}(x)$.

Now, for $y\in \OO(x)$, we consider the set
$$L_y := \{w\in \mathcal{O}(x) \mid w \leq y\},$$
which by the above paragraph and the transitivity of the action of $G$ on $\OO(x)$ must contain at least one element other than $y$.  Again by transitivity, $L_y$ cannot have a minimal element.  Therefore, $L_y$ must be infinite.  Additionally each pair of elements in $L_y - \{y\}$ have $y$ as an upper bound, so by the assumption that no pair of incomparable elements in $\OO(x)$ share an upper bound, $L_y$ must be totally ordered.  Therefore, the set
\begin{eqnarray}
\rho_y := \bigcup_{w \in L_y} GS_{(y,w)} \label{EQ:rhoy}
\end{eqnarray}
can be written as an infinite increasing union of geodesic spines, so it is a ray.

We claim that for any $s,t\in \OO(x)$ $\rho_s \cap \rho_t$ contains $\rho_w$ for some $w \in \OO(x)$.  To prove this, there are three cases to consider;  $s<t, t<s$ and $s\sim t$.  In the first case, $\rho_s \subset \rho_s\cap\rho_t$ and in the second case, $\rho_t \subset \rho_s\cap\rho_t$.  In the third case, $s\ncl t$ since we assumed that $x\ncl y$ for every incomparable pair $x\sim y$ in $\OO(x)$.  Since $\OO(x)$ is strongly connected, $s$ and $t$ share a lower bound in $\OO(x)$, say $w$.  In this case, $\rho_w \subset \rho_s \cap \rho_t$, proving the claim.

Let $R$ be the set of rays of the form given in Equation~(\ref{EQ:rhoy}).  We now prove that every ray in $R$ is infinite.  Note that, for any $g\in G$, $gL_y = L_{g\cdot y}$, so $g\rho_y = \rho_{g\cdot y}$, and  $G$ transitively permutes the set $R$.  Therefore, either all rays in $R$ are infinite or all are finite.  Suppose that they are all finite.  Since every two rays in $R$ eventually overlap, every two rays in $R$ have exactly the same set of endpoints.  Let $E$ be the set of endpoints of any (hence every) ray in $R$.  Then $G$ permutes the components of $T - E$.  Only one component contains points of $\OO(x)$, and $G$ must fix that component, which is therefore an invariant implicit subtree.  Since $G$ acts minimally on $T$, this is impossible, so the rays of $R$ must be infinite.

Since all rays in $R$ are infinite and any two eventually overlap, they all define the same end $\epsilon$, which is fixed by $G$.  Since we assumed $G$ not to fix a unique end of $T$, $G$ must fix another end $\delta$ of $T$.  By Lemma~\ref{L:twofixedends}, $G$ fixes the implicit line defined by $\epsilon$ and $\delta$, which by the minimality of the action, must be the entire manifold.  Therefore, $T \approx \mathbb R$, and $G$ acts faithfully on $\mathbb R$.  Thus, the order on $\OO(x)$ is a total order, and we have reached our desired contradiction.
\end{proof}

Now although we are unable to entirely omit restrictions on the stabilizer of the point $x$ we choose, we can weaken the restriction that the there be a point with trivial stabilizer.

\begin{thm} \label{T:nontrivialstab}
If a countable group $G$ acts on $T$ as in Theorem~\ref{T:groupgetsorder}, and there is some point $x \in T$ where $Stab(x)$ is left-orderable, then  $G$ admits a non-trivial left-invariant rectifiable simply connected partial order.
\end{thm}

\begin{proof}
By the same reasoning in Theorem~~\ref{T:groupgetsorder}, we see that the left cosets of $Stab(x)$ admit a rectifiable simply connected partial order. This extends to an order on the group $G$ as follows. Choose $g_1 \neq g_2 \in G$. If $g_1^{-1}g_2 \notin H$, then $g_1H \neq g_2H$, and we assign $g_1<g_2$ if $g_1H < g_2H$ in the partial order on the cosets, and similarly $g_2<g_1$ if $g_2H < g_1H$. Note that if such a pair has not been assigned to be comparable, then the cosets are also not comparable, so they have either an upper or a lower bound, which in turn provides $g_1$ and $g_2$ with common upper or lower bounds. If, on the other hand, If $g_1^{-1}g_2 \in H$, then either $e<g_1^{-1}g_2$ or $e >g_1^{-1}g_2$. In the first case we set $g_1 < g_2$ and in the second case we set $g_2 < g_1$. It is easy to see that the result is a left invariant rectifiable simply connected partial order. It cannot be trivial, since already at the level of cosets the order was not trivial.
\end{proof}

Next we move to the more general case of a minimal action on an oriented order tree. Again, we begin with actions which have at least one trivial stabilizer.

\begin{thm}\label{T:arbitrarytree}
If a countable group $G$ acts minimally and without fixing a unique end on an oriented order tree $T$ by orientation preserving homeomorphisms, and there is some point with trivial stabilizer, then $G$ admits a partial order with nontrivial left-invariant rectifiable simply connected extension.
\end{thm}

First, we have two lemmas necessary for the proof of Theorem~\ref{T:arbitrarytree}.  Recall the map, $\varphi: T' \to T$ mentioned after the proof of Proposition~\ref{P:tprime}.

\begin{lemma}\label{L:not.surj}
Let $T$ be an oriented order tree, let $T'$ be the associated $1$-manifold, and let $\hat T$ be the core of $T'$.  If $x \in \hat T$ and $X_1$ is a component of $T' - \{x\}$ then $\varphi(X_1) \neq T$.
\end{lemma}

\begin{proof}
Since $T'$ is a simply connected $1$-manifold, the point $x$ disconnects $T'$ into two components $X_1$ and $X_2$.  Let $y:= \varphi(x)$.  We first claim that $X_2 \not\subset \varphi^{-1}(y)$.  There are two possibilities for $y$.  Either $y$ is a branch point or a regular point of $T$.  If $y$ is regular, then $\varphi^{-1}(y)$ is the single point $x$, so  $X_2 \not\subset \varphi^{-1}(y)$.  If $y$ is a branch point, then the preimage of $y$ depends on the in-degree $n_o(y)$ and the out-degree $n_f(y)$.  We consider the case $n_o(y) = 0$ and $n_f(y) \geq 2$;  that is $y$ is a sink.  The other cases are similar. As in section 5 of \cite{RSS}, $n_f(y)= |R(y,f)| \geq 2$, where $R(y,f)$ is the set of incoming rays at $y$, where  such a ray  is an equivalence class of segments $\sigma$ with $f(\sigma)=y$ and where $\sigma_1 \approx \sigma_2 \Leftrightarrow \{y \} \subsetneqq \sigma_1 \cap \sigma_2$.  Now, $\varphi^{-1}(y)$ consists of an entire open ray $\sigma_y$ (the distinguished ray) and a set of points $\{x_{r_{\sigma}} \}$, where $r_{\sigma}$ ranges over all of the rays in $R(y,f)$. Since $n_f(y) \geq 2$, there are at least two of these points, one $x=x_{r_{\sigma}}$ and one $x_{r_{\tau}} \neq x$, where $\sigma$ and $\tau$ are segments representing $r_{\sigma}$ and $r_{\tau}$, respectively. Then $\sigma - \{ x \}$ and $\tau$ lie in different connected components of $T' - \{ x \}$, so $X_2$ contains either $\sigma - \{ x \}$ or $\tau$. But both $\sigma - \{ x \}$ and $\tau$ contain points not in $\varphi^{-1}(y)$, so in the case that $n_o(y) = 0$ and $n_f(y)\geq 2$, we have $X_2 \not\subset \varphi^{-1}(y)$.  The other possibilities for $n_o(y)$ and $n_f(y)$ are similar.

Now suppose towards a contradiction that $\varphi(X_1) = T$ and let $\alpha \in X_2$.  Then there exists $\beta \in X_1$ such that $\varphi(\beta) = \varphi(\alpha)$.  Since $\varphi^{-1}(\varphi(\alpha))$ is an implicit subtree of $T'$ and since $\alpha, \beta \in \varphi^{-1}(\alpha)$, we have $GS_{(\alpha,\beta)} \subset \varphi^{-1}(\varphi(\alpha))$.  Since $x \in GS_{(\alpha,\beta)}$, $\varphi(\alpha) = \varphi(x)=y$.  This is true for any point of $X_2$, so $X_2 \subset \varphi^{-1}(y)$, contradicting the previous claim.
\end{proof}

\begin{lemma}\label{L:t'action}
Suppose that $G$ acts minimally on the oriented order tree $T$ and let $T'$ be the associated $1$-manifold with the $G$-action.  Every nonempty invariant implicit subtree of $T'$ contains the core $\hat T$.
\end{lemma}
\begin{proof}
If $I$ is an invariant implicit subtree that does not contain $\hat T$, choose $x \in \hat T - I$.  Then $I$ is contained in one component $X_1$ of $T' - \{x\}$.  Since $\varphi(X_1) \neq T$, $\varphi(I) \neq T$.  But, $\varphi(I)$ is an invariant subtree of $T$, contradicting minimality of the action of $G$ on $T$.
\end{proof}

We now are in a position to prove Theorem~\ref{T:arbitrarytree}.

\begin{proof}[Proof of Theorem~\ref{T:arbitrarytree}]
Choose a point $\alpha \in T$ with trivial stabilizer. Blow $T$ up in the usual way to an oriented 1-manifold $T'$ on which $G$ acts. In the process, some points of $T$ may be split apart, and new rays and open intervals may be added. Then the points of $T'$ form a simply connected partially ordered set. Since $\varphi(\hat T) = T$, we may choose a point $x\in \hat T$ such that $\varphi(x)=\alpha$. Then $x$ will have trivial stabilizer as well.  As in Theorem~\ref{T:groupgetsorder}, we define a left-invariant partial order on $G$ by identifying it with $\OO(x)$. It may, of course, be the case that some incomparable pairs in $\OO(x)$ have no common bounds in $\OO(x)$, but they can be assigned the type $\sim_u$ or $\sim_l$ according to how they relate in $T'$. The resulting extension will clearly satisfy both Definitions~\ref{p:bounds} and~\ref{p:nobadtriads}, so it is a simply connected extension.  Since the poset $T'$ is rectifiable, the partial order of $\OO(x)$ is as well. The rest of the proof is devoted to showing that this extension is nontrivial.

We follow the proof of Theorem~\ref{T:quotientgetsorder}, using Lemma~\ref{L:t'action} instead of minimality to show that if the order is trivial, there must be two comparable elements of $\OO(x)$. Again, we then consider the set $R$ of rays of the form of equation~\ref{EQ:rhoy}, and note that $G$ acts transitively on $R$, so either all rays are infinite or all are finite.

The fact that, unlike in Theorem~\ref{T:groupgetsorder}, $G$ need not act minimally on $T'$ will complicate the rest of the proof.  Suppose that all the rays are finite, and let $E$ denote the set of endpoints of any (hence every) ray in $R$.  $E$ is $G$-invariant, and if no two points of $E$ are separable from each other, $E$ itself is an implicit subtree of $T'$.  In this case, $E$ certainly cannot contain all of $\hat T$, contradicting Lemma~\ref{L:t'action}.  If there are two separable points in $E$ then as in Theorem~\ref{T:groupgetsorder}, the set
$$I := \left(\bigcup_{a,b \in E}GS_{(a,b)} \right) - \mathcal{O}(x)$$
is a $G$-invariant implicit subtree of $T'$.  Since $\mathcal{O}(x) \subset \hat T$, $I$ does not contain $\hat T$, again contradicting Lemma~\ref{L:t'action}.  Therefore, all rays of $R$ must be infinite and as Theorem~\ref{T:groupgetsorder} they all define the same end $\epsilon'$, which is fixed by $G$.

We now use $\epsilon'$ to find an invariant implicit subtree or a fixed end of $T$.  First note that if $\omega$ is a geodesic spine in $T'$ then $\varphi(\omega)$ is a single point or a geodesic spine in $T$.  Since $\varphi$ maps points of $\mathcal{O}(x)$ to distinct points in $T$, $\varphi(\rho_y)$ is a ray for any $\rho_y \in R$.  Since any two rays in $R$ eventually overlap, the same is true of any two rays of the form $\varphi(\rho_y)$.  Moreover, $G$ transitively permutes the rays $\varphi(\rho_y)$.  So, either all are infinite or all are finite.  If they are finite, they all have the same set of endpoints, say $E_1$, which is a proper $G$-invariant implicit subtree of $T$.  Therefore, all rays $\varphi(\rho_y)$ are infinite, and they all define the same end $\epsilon$, which is fixed by $G$.

Since $G$ was assumed not to fix a unique end of $T$, there must be another fixed end $\delta$.  By Lemma~\ref{L:twofixedends}, $G$ fixes the implicit line $l$ defined by $\epsilon$ and $\delta$.  Since $G$ does act minimally on $T$, $l$ must be the entire tree, $T$.  Therefore, $T \approx \mathbb R$, and we have reached our desired contradiction.
\end{proof}

Just as in the case of a minimal action on a simply connected 1-manifold, the assumption on the stabilizer of $x$ can be weakened. The proof is essentially the same as for Theorem~\ref{T:nontrivialstab}.

\begin{thm}\label{T:weakstab}
If a countable group $G$ acts on $T$ as in Theorem~~\ref{T:arbitrarytree}, and there is some point $x \in T$ where $Stab(x)$ is left-orderable, then  $G$ admits a partial order with left-invariant nontrivial rectifiable simply connected extension.
\end{thm}



\begin{thebibliography}{999}

\bibitem{AB}
R. Alperin and H. Bass, \emph{Length functions of group actions on $\Lambda$-Trees}, Annals of Math. Studies, Princeton Univ. Press \textbf{111} (1987), 265-378.

\bibitem{Ba}
T. Barbot, \emph{Actions de groupes sur les 1-vari\'{e}t\'{e}s non s\'{e}par\'{e}es et feuilletages de codimension un} (French)[Actions of groups on non-Hausdorff 1-manifolds and codimension one foliations] Ann. Fac. Sci. Toulouse Math. \textbf{(6)7} (1998), no.4, 559-597.

\bibitem{bow}
B. Bowditch, \emph{Treelike structures arising from continua and convergence groups}, Mem. Amer. Math. Soc. \textbf{662} (1999), 1-86.

\bibitem{bc}
B. Bowditch and J. Crisp, \emph{Archimedean actions on median pretrees}, Math. Proc. Cambridge Philos. Soc. \textbf{130(3)} (2001), 383-400.

\bibitem{BRW}
S. Boyer, D. Rolfsen and B. Weist, \emph{Orderable $3$-manifold Groups}, arXiv:math.GT/0211110.

\bibitem{CM}
M. Culler and J.W. Morgan, \emph{Group actions on $\mathbb {R}$-trees}, Proc. London Math. Soc. \textbf{(3)55} (1987), 571-604.

\bibitem{dG94}
D. Gabai,\emph{Eight problems in the geometric theory of foliations and laminations on 3-Manifolds}, Georgia Conference Proceedings, Studies in Adv. Math. \textbf{2(2)} (1997), 1--33.

\bibitem{GK1}
D. Gabai and W. Kazez, \emph{Order trees and essential laminations of the plane},
Mathematical Research Letters \textbf{4} (1997), 603--616.

\bibitem{GO89}
D. Gabai and U. Oertel, \emph{Essential laminations in 3-manifolds}, Ann. Math. \textbf{130} (1989), 41--73.

\bibitem{Gh}
E. Ghys, \emph{Groups acting on the circle}, Ens. Math. \textbf{47} (2001), 329--407.

\bibitem{KM}
V. Kopytov and N. Medvedev, \emph{Right-Ordered Groups}, Siberian School of Algebra and Logic, Consultants Bureau, Plenum {1996}

\bibitem{p96}
F. Paulin, \emph{Actions de groupes sur les arbres}, S$\acute{e}$m. Bourbaki 1995-96, no. 808.

\bibitem{RS2}
R. Roberts and M. Stein, \emph{Group actions on order trees}, TopApp \textbf{115} (2001), 175--201.

\bibitem{RSS}
R. Roberts, J. Shareshian, and M. Stein, \emph{Infinitely many hyperbolic 3-manifolds which contain no reebless foliation}, J. A.M.S. \textbf{16(3)} (2003), 639--679.

\bibitem{sh1}
P. Shalen, \emph{Dendrology of groups: an introduction}, Essays in group theory (S.M. Gersten ed.), M.S.R.I. Pub. \textbf{8}, Springer-Verlag, 1987.



\end{thebibliography}
\end{document}